\documentclass[USenglish,journal]{IEEEtran}\pdfobjcompresslevel=0\pdfminorversion=4

\PassOptionsToPackage{table}{xcolor}
\PassOptionsToPackage{abbr,appendix}{clevethm}
\PassOptionsToPackage{compress}{cleveref}

\usepackage[
]{myPreamble}
\usepackage{cite}
\DeclareUnicodeCharacter{00A0}{ }
\setlist[itemize]{%
	label={\textbullet},
	itemsep=2pt,
	leftmargin=*,
	itemindent=\widthof{\textbullet\hspace*{\labelsep}},
}

	\newcommand\tf{{\vphantom f\smash{\hat f}}}%
	\newcommand\Ell{\mathcal L}%
	\newcommand\lspoint{\tilde}
	\newcommand\cont{\mathcal C}%
	\newcommand\geninv[1]{#1^\dag} 

	\newcommand\algname{Newton-type Alternating Minimization Algorithm}
	\newcommand\algacronym{NAMA}

	\usetikzlibrary{patterns}

	\pgfplotscreateplotcyclelist{mycolorlist}{%
		blue,  every mark/.append style={fill=blue!80!black},  mark=*\\%
		olive, every mark/.append style={fill=olive!80!black},   mark=square*\\%
		black, every mark/.append style={fill=black}, mark=triangle*\\%
		red,  every mark/.append style={fill=white},  mark=*\\%
		orange,  every mark/.append style={fill=white},  mark=square*\\%
		red!50!blue,   every mark/.append style={fill=white},  mark=triangle*\\%
		gray, every mark/.append style={fill=white},  mark=diamond*\\%
	}

	\pgfplotsset{every tick label/.append style={font=\footnotesize}}
	\pgfplotsset{every axis label/.append style={font=\footnotesize}}
	\pgfplotsset{legend style={font=\footnotesize}}

\let\relint\relax
\DeclareMathOperator\relint{ri}

\begin{document}

	\title[Newton-type Alternating Minimization Algorithm]{%
		Newton-type Alternating Minimization Algorithm for Convex Optimization%
	}
	\author{%
		Lorenzo Stella,
		Andreas Themelis
		and
		Panagiotis Patrinos
		\thanks{%
			All authors are affiliated with the \TheAddressKU.
			The first two authors are also affiliated with the \TheAddressIMT.\newline
			\{%
				\href{mailto:lorenzo.stella@esat.kuleuven.be}{lorenzo.stella},
				\href{mailto:andreas.themelis@esat.kuleuven.be}{andreas.themelis},
				\href{mailto:panos.patrinos@esat.kuleuven.be}{panos.patrinos}%
			\}%
			\href{mailto:lorenzo.stella@esat.kuleuven.be}{@esat.kuleuven.be}.
		}%
	}%

	\maketitle
	\begin{abstract}
		We propose NAMA (Newton-type Alternating Minimization Algorithm) for solving structured nonsmooth convex optimization problems where the sum of two functions is to be minimized, one being strongly convex and the other composed with a linear mapping.
The proposed algorithm is a line-search method over a continuous, real-valued, exact penalty function for the corresponding dual problem, which is computed by evaluating the augmented Lagrangian at the primal points obtained by alternating minimizations.
As a consequence, NAMA relies on exactly the same computations as the classical alternating minimization algorithm (AMA), also known as the dual proximal gradient method.
Under standard assumptions the proposed algorithm possesses strong convergence properties,
while
under mild additional assumptions the asymptotic convergence is superlinear, provided that the search directions are chosen according to quasi-Newton formulas.
Due to its simplicity, the proposed method is well suited for embedded applications and large-scale problems.
Experiments show that using limited-memory directions in NAMA greatly improves the convergence speed over AMA and its accelerated variant.

	\end{abstract}


	\section{Introduction}
		\label{sec:Introduction}
		We consider convex optimization problems of the form
\[\tag{P}\label{eq:Problem}
	\minimize_{x\in\R^n}\ f(x) + g(Ax),
\]
where \(f\) is strongly convex, \(g\) is convex and \(A\) is a linear mapping.
Problems of this form are quite general and appear in various areas of applications, including optimal control \cite{stathopoulos2016splitting}, system identification \cite{fazel2013hankel} and machine learning \cite{boyd2011distributed,parikh2014proximal}. For example, whenever \(g\) is the indicator function of a convex set \(C\), then \eqref{eq:Problem} models a constrained convex problem: if \(C\) is a box, then in particular \eqref{eq:Problem} amounts to minimizing a strongly convex function subject to polyhedral constraints.

A general approach to the solution of \eqref{eq:Problem} is based on the dual proximal gradient method, or forward-backward splitting, also known as alternating minimization algorithm (AMA) \cite{tseng1991applications}.
This is the dual application of an algorithm introduced by Lions and Mercier \cite{lions1979splitting} for finding the zero of the sum of two maximal monotone operators, one of which is assumed to be co-coercive.
The alternating minimization algorithm is intimately tied to the framework of \emph{augmented Lagrangian} methods, and its global convergence and complexity bounds are well covered in the literature, see \cite{tseng1991applications}:
a global convergence rate of order \(O(1/\sqrt{k})\) holds for the primal iterates of AMA under very general assumptions, and can be improved to the optimal rate \(O(1/k)\) using a simple acceleration technique due to Nesterov, see \cite{beck2009fast,nesterov2013gradient,beck2014fast}.

As with all first order methods, the performance of (fast) AMA is severely affected by ill-conditioning of the problem \cite{stathopoulos2016splitting}.
One way to deal with this issue, which is extensively used in classical smooth, unconstrained optimization, is to precondition the problem using (approximate) second-order information on the cost function, as in (quasi-) Newton methods.
However, both \eqref{eq:Problem} and its dual are nonsmooth in general.
This motivates considering the concept of \emph{alternating minimization envelope} (AME): this is a real-valued (as opposed to \emph{extended} real-valued) exact merit function for the dual problem, and is precisely the augmented Lagrangian associated with \eqref{eq:Problem} evaluated at the primal points computed by AMA.
Under mild assumptions on \eqref{eq:Problem}, the AME is continuously differentiable around the set of dual solutions and even strictly twice differentiable there. As a consequence, the AME allows to extend classical, smooth unconstrained optimization algorithms to the solution of the dual problem to \eqref{eq:Problem}, which is nonsmooth in general.
In this work we propose a dual line-search method, which uses the AME as merit function to compute the stepsizes. The convergence properties of the proposed algorithm greatly improve over AMA when fast-converging directions, computed by means of quasi-Newton formulas, are followed.
Furthermore, we show that the AME is equivalent to the \emph{forward-backward envelope} (FBE, see \cite{patrinos2013proximal, stella2017forward, themelis2016forward}) associated with the dual problem.

\subsection{Related works}

The FBE, as a tool for extending smooth unconstrained algorithms to nonsmooth problems, has first been introduced in \cite{patrinos2013proximal}: there, two semismooth Newton methods are proposed for minimizing the sum of two convex functions, one of which is smooth and the other having an efficiently computable proximal mapping. This is the classical setting in which the proximal gradient method (and its accelerated variant) can be applied.
In \cite{stella2017forward} the convexity assumption on the smooth term is relaxed, and the authors propose a line-search method with
global sublinear rate (in the convex case) and asymptotic superlinear rate when quasi-Newton directions are used: the algorithm
relies on descent directions over the FBE which is required to be \emph{everywhere} differentiable. In \cite{liu2016further} classical gradient-based line-search methods are considered for minimizing the FBE, see also \cite{sampathirao2017proximal}.
In \cite{themelis2016forward} the most general framework, where both summands are allowed to be nonconvex, is taken into account. In this case differentiability of the FBE cannot be assumed: a new algorithm is proposed which computes fast convergent directions with no need for gradient information on the FBE.

A similar approach was used in \cite{patrinos2014douglas,themelis2017douglas} to accelerate other splitting algorithms, namely the \emph{Douglas-Rachford splitting} and its dual counterpart ADMM.

\subsection{Contributions and organization of the paper}

In the present paper we deal with the case where \(g\) in \eqref{eq:Problem} is composed with a linear mapping. In this case, even though \(g\) may possess an efficiently computable proximal mapping, \(g\circ A\) in general does not. This motivates addressing the dual problem of \eqref{eq:Problem} instead.
The contributions and organization of the present work can be summarized as follows.

\begin{itemize}
	\item We propose the \DEF{\algname} (\algacronym, \Cref{sec:Background}, \Cref{alg:NAMA}), a generalization of the alternating minimization algorithm that performs a line-search step over the AME: the proposed algorithm relies on the very same alternating minimization operations of AMA.
	\item We show that the AME is equivalent to the FBE of the dual problem (\Cref{sec:DualFBE}). This observation extends a classical result by Rockafellar, relating the Moreau envelope and the augmented Lagrangian, to our setting where an additional strongly convex term is present.
	\item We show that the proposed method enjoys global sublinear convergence under standard assumptions, and local linear convergence assuming \emph{calmness} of the subdifferentials of the problem terms (\Cref{sec:Convergence}).
	\item We analyze the first- and second-order properties of the AME, by linking them to generalized second-order properties of the primal functions \(f\) and \(g\) (\Cref{sec:DiffFBE}).
	\item We show that the proposed method converges asymptotically superlinearly when the dual problem has a (unique) strong dual minimum, and the line-search directions are selected so as to satisfy the Dennis-Mor\'e condition, as it is the case when quasi-Newton update formulas are adopted (\Cref{sec:Superlinear}). The effectiveness of our approach is demonstrated by numerical simulations on linear MPC problems (\Cref{sec:Simulations}).
\end{itemize}

Differently from the approaches in \cite{stella2017forward,liu2016further,sampathirao2017proximal}, \algacronym{} does not require the gradient of the envelope function, therefore no second-order information on the smooth term is needed: this would severely limit its applicability in the present setting where the dual problem is solved.
Furthermore, with respect to the approaches of \cite{liu2016further,sampathirao2017proximal}, the algorithm presented here possesses strong global convergence properties which are not typical of classical line-search methods.
Differently from \cite{themelis2016forward}, despite the fact that the selected directions may not be descent directions and the line search is performed on the envelope function, \algacronym{} is a descent method for the dual objective: this allows to simplify the convergence analysis of the method, and to show the global sublinear convergence rate for the dual cost and the primal iterates.


		\subsection{Notation}
			\label{sec:Notation}
			In what follows \(\innprod{\cdot}{\cdot}\) denotes an inner product over a Euclidean space (whose nature will be clear from the context) and \(\|\cdot\| = \sqrt{\innprod{\cdot}{\cdot}}\) is the associated norm.
For a linear \(A:\R^n\to\R^m\), \(\|A\|\) is the operator norm induced by the inner products over \(\R^n\) and \(\R^m\).
For a set \(C\), we denote by \(\relint(C)\) its relative interior, and by \(\proj_C(x) = \argmin_{y\in C}\|y-x\|\) the projection onto \(C\) in the considered norm.
We denote the extended real line by \(\overline\R = \R\cup\set{\infty}\),
and by \(\Gamma_0(\R^n)\) the set of proper, closed, convex functions defined over \(\R^n\) with values in \(\Rinf\).
For \(h\in\Gamma_0(\R^n)\) its \DEF{Fenchel conjugate} \(\conj{h}\), defined as
	\(\conj h(y)
	=
	{\sup}_{x\in\R^n}
	\set{\innprod{x}{y} - h(x)}\) is also proper, closed and convex.
Properties of conjugate functions are well described for example in~\cite{rockafellar1997convex,hiriart2001fundamentals,bauschke2011convex,rockafellar2011variational}.
Among these we recall the Fenchel-Young inequality \cite[Prop. 13.13]{bauschke2011convex}
\begin{equation}\label{eq:FenchelIneq}
	\innprod{x}{y} \leq h(x) + \conj h(y)\quad \forall x,y,
\end{equation}
with
\begin{equation}\label{eq:ConjSubgr}
	y\in\partial h(x) \Leftrightarrow \innprod{x}{y} = h(x) + \conj h(y) \Leftrightarrow x\in\partial \conj h(y),
\end{equation}
see \cite[Thm. 23.5]{rockafellar1997convex}.
For any \(\gamma>0\), the \DEF{proximal mapping} associated with $h$, with stepsize \(\gamma\), is denoted as
\[
	\prox_{\gamma h}(x) = {\argmin}_{z}\set{h(z)+(\nicefrac{1}{2\gamma})\|z-x\|^2}.
\]
This satisfies the Moreau identity \cite[Thm. 14.3(ii)]{bauschke2011convex}
\begin{equation}\label{eq:MoreauId}
	y = \prox_{\gamma h}(y) + \gamma \prox_{\gamma^{-1} \conj h}(\gamma^{-1}y)
		\quad \forall y.
\end{equation}
The value function of the problem defining \(\prox_{\gamma h}\) is the \emph{Moreau envelope}
\[ h^\gamma(x) = {\min}_{z}\set{h(z)+(\nicefrac{1}{2\gamma})\|z-x\|^2}. \]
An alternative formulation for \eqref{eq:Problem} is
\begin{equation}\label{eq:AltProblem}\tag{P$^\prime$}
	\minimize_{x\in\R^n,z\in\R^m}f(x) + g(z)
\quad
	\stt Ax = z.
\end{equation}
Therefore we can define the \DEF{augmented Lagrangian} associated with \eqref{eq:Problem}, denoted as
\[
	\Ell_\gamma(x,z,y)
{}={}
	f(x) + g(z) + \innprod{y}{Ax-z} + \tfrac{\gamma}{2}\|Ax-z\|^2,
\]
where \(\gamma\geq 0\).
We indicate by \(\Ell \equiv \Ell_0\) the ordinary Lagrangian function.


We follow the terminology of \cite{rockafellar2011variational} when referring to the concepts of \emph{strict continuity} and \emph{strict differentiability}. We say that a mapping \(F:\R^n\to\R^m\) is \emph{strictly continuous} at \(\bar x\) if \cite[Def. 9.1(b)]{rockafellar2011variational}
\[
	\limsup_{
		\substack{
			(x,y)\to(\bar x,\bar x)\\
			x\neq y
		}
	}{
		\frac{\smash{\|F(y)-F(x)\|}}{\|y-x\|}
	}
{}<{}
	\infty.
\]
If \(F\) is (Frech\'et) differentiable, we let \(\jac F{}:\R^n\to\R^{m\times n}\) denote the \emph{Jacobian} of \(F\).
When \(m=1\) we indicate with \(\nabla F=\jac F{}^\top\) the \emph{gradient} of \(F\) and with \(\nabla^2 F=\jac{\nabla F}{}^\top\) its \emph{Hessian}, whenever it makes sense.
We say that $F$ is \emph{strictly differentiable} at $\bar x$ if it satisfies the stronger limit \cite[Eq. 9(7)]{rockafellar2011variational}
\[
	\lim_{
		\substack{
			(x,y)\to(\bar x,\bar x)\\
			x\neq y
		}
	}{
		\frac{\|F(y)-F(x)-\jac F{\bar x}[y-x]\|}{\|y-x\|}
	}
{}={}
	0.
\]
Some results in the paper are based on generalized second-order properties of extended-real-valued functions.
\begin{defin}[\hspace*{-3pt}{\cite[Def. 13.6]{rockafellar2011variational}}]%
Function \(h : \R^n\to\Rinf\) is said to be \emph{twice epi-differentiable} at \(x\) for \(v\), if the second-order difference quotient 
\[
	\twicequot[v]{h}[\tau]{x}[d]
{}={}
	\frac{
		\smash{
			h(x+\tau d)-h(x)-\tau\innprod{v}{d}
		}
	}{
		\nicefrac{\tau^2}{2}
	}
\]
epi-converges as \(\tau \searrow 0\) (\ie its epigraph converges in the sense of Painlev\'e-Kuratowksi, see~\cite[Def.~7.1]{rockafellar2011variational}),
the limit being the function \(\twiceepi[v]{h}{x}\) given by
\begin{equation*}\label{eq:INTRO:2ndSubder}
	\twiceepi[v]{h}{x}[d]
=
	\liminf_{\substack{\tau\searrow 0\\d'\to d}} \twicequot[v]{h}[\tau]{x}[d'].
\end{equation*}
In this case \(\twiceepi[v]{h}{x}[d]\), as a function of \(d\), is said to be the \emph{second-order epi-derivative} of $h$ at $x$ for $v$.
If \(\twicequot[\bar v]{h}[\tau]{\bar x}\) epi-converges as \(\tau \searrow 0\), \(\bar x \to x\) and \(\bar v \to v\), then \(h\) is said to be \emph{strictly} twice epi-differentiable.
\end{defin}
Twice epi-differentiability is a mild requirement, and functions with this property are abundant.
Refer to \cite{rockafellar1988first, rockafellar1989second, poliquin1992amenable, poliquin1995second, poliquin1996generalized}
and to \cite[§7, §13]{rockafellar2011variational} for examples and an in-depth account on epi-derivatives, epi-differentiability, and their connections with ordinary differentiability.


	\section{Background and proposed algorithm}
		\label{sec:Background}
Without further specifying it, throughout the paper we will work under the following basic assumption.
\begin{ass}\label{ass:fg_TEMP}\label{ass:Basic}
The following hold for \eqref{eq:Problem}:
\begin{enumerate}
	\item\label{ass:feasibility}
		\eqref{eq:Problem} is feasible, \ie \(A\dom f \cap \dom g\neq \emptyset\);
	\item\label{ass:fStrConvex}
		\(f\in\Gamma_0(\R^n)\) is strongly convex with modulus \(\mu_f>0\);%
		\footnote{Function \(h\) has convexity modulus \(c\geq 0\) if \(h-\frac c2\|{}\cdot{}\|^2\) is convex.}
	\item\label{ass:gConvex}
		\(g\in\Gamma_0(\R^m)\).
\end{enumerate}
\end{ass}
\begin{rem}\label{rem:Basic}
\cref{ass:Basic} guarantees, by strong convexity of \(f\), that a solution to \eqref{eq:Problem} exists and is unique, be it \(x_\star\).
\cref{ass:fStrConvex} also implies that \(\conj f\) is Lipschitz continuously differentiable with constant \(\mu_f^{-1}\) \cite[Th. 12.60]{rockafellar2011variational}.
\cref{ass:gConvex} ensures that \(\conj g\) is also proper, closed, convex \cite[Cor. 13.33]{bauschke2011convex},
and its Moreau envelope \((\conj g)^\gamma\) is strictly continuous \cite[Ex. 10.32]{rockafellar2011variational} with \(\gamma^{-1}\)-Lipschitz gradient
\begin{equation}\label{eq:MoreauGradient}
	\nabla(\conj g)^\gamma(y)
{}={}
	\gamma^{-1}\left(y-\prox_{\gamma\conj g}(y)\right),
\end{equation}
as shown in \cite[Prop. 12.29]{bauschke2011convex}.
\end{rem}

The Fenchel dual problem associated with \eqref{eq:Problem} is
\[\tag{D}\label{eq:DualProblem}
	\minimize_{y\in\R^m}\ \psi(y)
{}={}
	\conj f(-\trans{A} y)
	{}+{}
	\conj g(y).
\]
Under \cref{ass:Basic} strong duality holds, see \cite[Thm. 5.2.1(b)-(c)]{auslender2003asymptotic}
and primal-dual solutions \((x_\star,y_\star)\) to \eqref{eq:Problem}-\eqref{eq:DualProblem} are characterized by the first-order optimality conditions
\begin{subequations}\label{eq:FirstOrderNecessary}
\begin{align}
	-\trans A y_\star
{}\in{}
	\partial f(x_\star)
\quad &
	(%
		\Leftrightarrow\
			x_\star=\nabla\conj f(-\trans{A} y_\star)
	)
\\
	y_\star
{}\in{}
	\partial g(A x_\star)
\quad &
	(%
		\Leftrightarrow\
			A x_\star\in\partial\conj g(y_\star)
	).
\end{align}
\end{subequations}

A natural way to tackle \eqref{eq:Problem} is to solve \eqref{eq:DualProblem} by means of forward-backward splitting (or proximal gradient method): starting from an initial dual point \(y^0\in\R^m\), iterate
\begin{equation}\label{eq:FBS}
	y^{k+1}=T_{\gamma}(y^k) \coloneqq \prox_{\gamma \conj g}(y^k+\gamma A\nabla \conj f(-\trans{A}y^k))
\end{equation}
for some positive stepsize parameter \(\gamma\). If we define
the associated \DEF{fixed-point residual}%
\[
	R_\gamma(y)
{}\coloneqq{} 
	\gamma^{-1}(y-T_\gamma(y)),
\]
then dual optimality can be characterized as follows:
\begin{equation}\label{eq:criticalOptimal}
	y_\star\in Y_\star
\;\Leftrightarrow\;
	y_\star\in\fix T_\gamma \;\Leftrightarrow\; y_\star\in\zer R_\gamma ~~\forall\gamma>0.
\end{equation}

Iterations \eqref{eq:FBS} are easily shown to be equivalent to the following scheme, the \emph{alternating minimization algorithm} (AMA)
\begin{subequations}\label{eq:AMM}
\begin{align}
\label{eq:AMMx}
	x^k
{}={} &
	\fillwidthof[c]%
	{z_{\gamma}(y^k)}%
	{x(y^k)}
{}={}
	\argmin_{x\in\R^n}\set{f(x) + \innprod{y^k}{Ax}},
&
\\
\label{eq:AMMz}
 	z^k
{}={} &
	z_{\gamma}(y^k)
{}={}
	\argmin_{z\in\R^m} \Ell_\gamma(x^k, z, y^k),
\\
\label{eq:AMMdual}
	y^{k+1}
{}={} &
	y^k + \gamma(Ax^k - z^k).
\end{align}
\end{subequations}
Note that step \eqref{eq:AMMz} can be equivalently formulated as
\[ z^k = \prox_{\gamma^{-1} g}(\gamma^{-1}y^k+Ax(y^k)). \]
Using the notation of \eqref{eq:AMM}, \(T_\gamma\) and \(R_\gamma\) can be expressed as
\ifieee
	\begin{subequations}\label{eq:TR}\begin{align}
		T_\gamma(y)
	{}={}&
		y + \gamma(A x(y) - z_\gamma(y)) \\
		R_\gamma(y)
	{}={}&
		z_\gamma(y) - A x(y).
	\end{align}\end{subequations}
\fi
\ifarxiv
	\begin{equation}\label{eq:TR}
		T_\gamma(y)
	{}={}
		y + \gamma(A x(y) - z_\gamma(y))
	\qquad\text{and}\qquad
		R_\gamma(y)
	{}={}
		z_\gamma(y) - A x(y).
	\end{equation}
\fi
It can be shown that \(x^k\to x_\star\) in iterations \eqref{eq:AMM}, provided that \(\gamma\in(0,\nicefrac{2\mu_f}{\|A\|^2})\), see \cite[Prop. 3]{tseng1991applications}.
Moreover, the dual cost in this case converges sublinearly to the optimum with global rate \(O(1/k)\), and the extrapolation techniques introduced by Nesterov~\cite{nesterov1983method,nesterov2003introductory,nesterov2013gradient} allow to obtain accelerated versions of AMA with an optimal global rate \(O(1/k^2)\), see \cite{beck2014fast}: we will here refer to this variant as \emph{fast} AMA.

		\subsection{Newton-type alternating minimization algorithm}
			\label{sec:NAMA}

\begin{algorithm}[t]
	\algcaption{Newton-type AMA (\algacronym)}%
	\label{alg:NAMA}%
	\begin{algorithmic}[1]
	\Require{
		\(y^0\in\R^m\),
		\(\gamma \in (0,\nicefrac{\mu_f}{\|A\|^2}) \),
		\(\beta\in(0,1)\)%
	}%
	\Initialize{%
		\(k = 0\)%
	}%
	%
	\STATE\label{step:AMLS1_altmin}%
		\(x^k = \argmin_{x} \set{f(x) + \innprod{y^k}{Ax}}\)%
	\item[]%
		\(\fillwidthof[r]{x^k}{z^k} = \argmin_{z} \Ell_\gamma(x^k, z, y^k)\)
	%
	\STATE\label{step:AMLS1_direction}%
		Choose a direction \(d^k\in\R^m\)%
	\STATE\label{step:AMLS1_linesearch}
		Find the largest \(\tau_k=\beta^{i_k}, i_k\in\N,\) such that
		\begin{equation}\label{eq:AMLS_ascent}
			\Ell_{\gamma}(\lspoint x^k, \lspoint z^k, \lspoint y^k) \geq \Ell_{\gamma}(x^k, z^k, y^k),
		\end{equation}
		where
		\begin{equation*}
			\begin{aligned}
				\lspoint y^k &{}={} y^k + \tau_k d^k + \gamma (1-\tau_k) (Ax^k-z^k)\\
				\lspoint x^k &{}={} {\argmin}_{x} \set{f(x) + \innprod{\tilde y^k}{Ax}}\\
				\lspoint z^k &{}={} {\argmin}_{z} \Ell_\gamma(\lspoint x^k, z, \lspoint y^k)
			\end{aligned}
		\end{equation*}
%
	%
	\STATE\label{step:AMLS1_global}
		\(y^{k+1} = \lspoint y^k + \gamma (A\lspoint x^k - \lspoint z^k)\),
		\(k = k+1\),
		go to step \ref{step:AMLS1_altmin}
\end{algorithmic}

\end{algorithm}

The convergence speed of (fast) AMA is affected by ill-conditioning of the problem, as it is the case for all first-order methods.
To accelerate convergence, we propose \Cref{alg:NAMA}.
An overview of the algorithm is as follows:
\begin{itemize}
	\item \Cref{alg:NAMA} is composed by the very same operations as AMA: in fact, only alternating minimization steps with respect to \(x\) and \(z\) are performed.
	\item Step \ref{step:AMLS1_linesearch} computes a new dual iterate \(\tilde{y}^k\), by performing a line-search over the augmented Lagrangian associated with \eqref{eq:Problem} evaluated at the alternating minimization primal points: we will see that this is equivalent to the forward-backward envelope function associated with the dual problem \eqref{eq:DualProblem}.
	\item The line-search is performed using a convex combination of the ``nominal'' residual direction \(\gamma(Ax^k-z^k)\) and an ``arbitrary'' direction \(d^k\), to be selected so as to ensure fast asymptotic convergence. This novel choice of direction ensures that the line-search is feasible at every iteration (\ie, condition \eqref{eq:AMLS_ascent} holds for a sufficiently small stepsize) despite the fact that \(d^k\) may not be a direction of descent, as we will see.
	\item Step \ref{step:AMLS1_global} will allow us to obtain global convergence rates, and it comes at no cost since vectors \(\tilde y^k, \tilde x^k, \tilde z^k\) have already been computed in the line-search. In a sense, this step robustifies the algorithmic scheme.
\end{itemize}

By appropriately choosing \(d^k\), the algorithm is able to greatly improve the convergence of AMA: we will prove that the algorithm converges with superlinear asymptotic rate when Newton-type directions are selected.
For this reason we refer to \Cref{alg:NAMA} as \emph{\algname{} (NAMA)}.


\begin{rem}[AMA as special case]\label{rem:GeneralizesAMA}
If in \Cref{alg:NAMA} one sets $d^k = 0$ for all $k$, then one can trivially select $\tau_k = 1$.
In this case,
\(
	(\tilde y^k, \tilde x^k, \tilde z^k)
{}={}
	(y^k, x^k, z^k)
\)
and \Cref{alg:NAMA} reduces to AMA, cf. \eqref{eq:AMM}.
\end{rem}

\begin{rem}[General equality constrained problems]\label{rem:GeneralECCP}
For any proper, closed, convex \(h:\R^r\to\Rinf\), \(b\in\R^m\)
and linear mapping \(\func{B}{\R^r}{\R^m}\), a problem of the form
\begin{equation}\tag{P$^{\prime\prime}$}\label{eq:OtherProblem}
	\minimize_{x\in\R^n,w\in\R^r} f(x) + h(w)
\quad
	\stt Ax+Bw=b
\end{equation}
can be rewritten as \eqref{eq:Problem} by letting
\begin{equation}\label{eq:substitution}
	g(z) = (Bh)(b-z) = \inf_{w\in\R^r}\set{h(w)}[Bw = b - z].
\end{equation}
Function \((Bh)\) is the \emph{image of \(h\) under \(B\)}, see \cite[Thm. 5.7]{rockafellar1997convex} and discussion thereafter.
If we further assume \(\relint(\dom \conj h)\cap\range(\trans{B}) \neq \emptyset\), then \((Bh)\) is proper, closed, convex, see
\cite[Thm. 16.3]{rockafellar1997convex},
therefore \(g\) in \eqref{eq:substitution} satisfies \cref{ass:gConvex} (if \(h\) is piecewise linear-quadratic then it is sufficient to assume \(\dom \conj h\cap\range(\trans{B}) \neq \emptyset\), see \cite[Cor. 11.33(b)]{rockafellar2011variational}).
In this case steps \eqref{eq:AMMz} and \eqref{eq:AMMdual} of AMA become
\begin{align*}
	 	w^k
	{}={} &
		\argmin_{w\in\R^r}\set{g(w) + \innprod{y^k}{Bw} + \tfrac{\gamma}{2}\|Ax^k + Bw - b\|^2}
	\\
		y^{k+1}
	{}={} &
		y^k + \gamma(Ax^k + Bw^k - b).
\end{align*}
Similar modifications allow to adapt NAMA to this more general setting:
in light of these observations, what follows readily applies to problems of the form \eqref{eq:OtherProblem}.
\end{rem}

		\subsection{Quasi-Newton directions}
			\label{sec:QN}
			There is freedom in selecting \(d^k\) in \Cref{alg:NAMA}.
To accelerate the convergence of the iterates, one possible choice is to compute fast converging directions for the system of nonlinear equations \(R_\gamma(y) = 0\) characterizing dual optimal points, cf. \eqref{eq:criticalOptimal}. Specifically, in \Cref{alg:NAMA} one can set
\begin{equation}\label{eq:QNdir}
	d^k = B_k^{-1}(Ax^k - z^k),
\end{equation}
for a sequence of nonsingular matrices \(\seq{B_k}\) approximating in some sense the Jacobian \(\jac{R_\gamma}{}\) at the limit point of the dual iterates \(\seq{y^k}\).
In quasi-Newton methods, starting from an initial nonsingular matrix \(B_0\), the sequence of matrices \(\seq{B_k}\) is determined by low-rank \emph{updates} that satisfy the secant condition: in \Cref{alg:NAMA} fast asymptotic convergence can be proved if
\[
	B_{k+1} p^k = q^k\quad\mbox{with}\quad
	\begin{cases}
		p^k {}={} & \lspoint y^k-y^k, \\
		q^k {}={} & (\lspoint z^k - A\lspoint x^k) - (z^k - Ax^k),
	\end{cases}
\]
as will be discussed in \Cref{sec:Superlinear}.
Note that all quantities required to compute the vectors \(p^k, q^k\) are available as by-product of the iterations.

In \cite{powell1970hybrid} the modified Broyden update is proposed, that prescribes rank-one updates of the form
\def\BFGSformula{B_k + \frac{q^k\trans{(q^k)}}{\innprod{q^k}{p^k}} - \frac{B_{k}p^k\trans{(B_{k}p^k)}}{\innprod{p^k}{B_{k}p^k}}.}%
\begin{equation}\label{eq:ModifiedBroyden}
	\text{\it
		Broyden%
	}
\qquad
	B_{k+1}
{}={}
	B_k
	{}+{}
	\theta_k\frac{(q^k-B_k p^k)\trans{(p^k)}}{\|p^k\|^2}.
\end{equation}
Here, \(\seq{\theta_k}\subset [0,2]\) is a sequence used to ensure that all terms in \(\seq{B_k}\) are nonsingular, so that \eqref{eq:QNdir} is well defined.
The original Broyden method \cite{broyden1965class} is obtained with \(\theta_k \equiv 1\).

Probably the most popular quasi-Newton scheme is BFGS, which prescribes the following rank-two updates
\begin{equation}\label{eq:BFGS}
	\text{\it BFGS}
\qquad
	B_{k+1}
{}={}
	\BFGSformula
\end{equation}
Note that in this case matrices \(B_k\) are symmetric, and in fact the fast asymptotic properties of BFGS are guaranteed only if the Jacobian \(\jac{R_\gamma}{}\) is symmetric \cite{byrd1989tool} at the problem solution.
This is not the case in our setting (cf. \Cref{ex:QP}) although we have observed that \eqref{eq:BFGS} often outperforms other non-symmetric updates such as \eqref{eq:ModifiedBroyden} in practice.

Using the Sherman-Morrison-Woodbury identity in \eqref{eq:ModifiedBroyden} and \eqref{eq:BFGS} allows to directly store and update \(H_k=B_k^{-1}\), so that \(d^k\) can be computed without inverting matrices or solving linear systems.

Ultimately, instead of storing and operating on dense \(m\times m\) matrices, \emph{limited-memory} variants of quasi-Newton schemes keep in memory only a few (usually \(3\) to \(30\)) most recent pairs $(p^k,q^k)$ implicitly representing the approximate inverse Jacobian.
Their employment considerably reduces storage and computations over the full-memory counterparts, and as such they are the methods of choice for large-scale problems. The most popular limited-memory method is probably L-BFGS, which is based on the update \eqref{eq:BFGS}, but efficiently computes matrix-vector products with the approximate inverse Jacobian using a \emph{two-loop recursion} procedure \cite{liu1989limited,nocedal1980updating,nocedal2006numerical}.

	\section{Alternating minimization envelope}
		\label{sec:DualFBE}
		The fundamental tool enabling fast convergence of \Cref{alg:NAMA} is the \emph{alternating minimization envelope} function associated with \eqref{eq:Problem}. This is precisely the (negative)
augmented Lagrangian function, evaluated at the primal points given by the alternating minimization steps.

\begin{defin}[Alternating minimization envelope]\label{def:AME}
The \DEF{alternating minimization envelope} (AME) for \eqref{eq:Problem}, with parameter \(\gamma>0\), is the function
\[
	\psi_\gamma(y) = -\Ell_\gamma(x(y),z_\gamma(y),y).
\]
\end{defin}

The first observation that we make relates the alternating minimization envelope in \Cref{def:AME} with the concept of \emph{forward-backward envelope}.

\begin{thm}\label{thm:AMEvsFBE}
	Function \(\psi_\gamma\) is the \DEF{forward-backward envelope} (cf. \cite[Def. 2.1]{stella2017forward}) associated with the dual problem \eqref{eq:DualProblem}:
	\begin{align}
		\nonumber
			\psi_\gamma(y)
		&{}={}
			\conj{f}(-\trans{A} y) + \conj{g}(T_\gamma(y)) + \tfrac{\gamma}{2}\|Ax(y)-z_\gamma(y)\|^2
		\\
		\label{eq:AMEvsFBE}
		&\phantom{{}={}}
			+\gamma\innprod{ Ax(y)}{z_\gamma(y) - Ax(y)}.
	\end{align}
	\begin{proof}
		The optimality conditions for \(x(y)\) and \(z_\gamma(y)\) are
	\begin{subequations}\label{subeq:OptCondxy}
		\begin{align}
			\label{eq:OptCondx}
			\partial f(x(y))
		&{}\ni{}
			-\trans{A} y,
		\\
			\label{eq:OptCondy}
			\partial g(z_\gamma(y))
		&{}\ni{}
			T_\gamma(y) = y+\gamma(Ax(y)-z_\gamma(y)).
		\end{align}
	\end{subequations}
	From these, using \eqref{eq:ConjSubgr}, we obtain
	\begin{subequations}\label{eq:ConjSubgr_fg}\begin{align}
		\label{eq:ConjSubgr_f}
			f(x(y))+\conj{f}(-\trans{A} y)
		&{}={}
			-\innprod{ Ax(y)}{y}
		\\
		\label{eq:ConjSubgr_g}
			g(z_\gamma(y))+\conj{g}(T_\gamma(y))
		&{}={}
			\innprod{z_\gamma(y)}{T_\gamma(y)}
	\end{align}\end{subequations}
	Summing \eqref{eq:ConjSubgr_fg} and rearranging the terms we get \eqref{eq:AMEvsFBE}.
\end{proof}
\end{thm}

An alternative expression for \(\psi_\gamma\) in terms of the Moreau envelope of \(\conj g\) is as follows, see \cite{patrinos2013proximal}:
\begin{equation}\label{eq:EquivFBE}
		\psi_\gamma(y)
	{}={}
		\conj f(-\trans Ay)
	{}-{}
		\tfrac{\gamma}{2}\|Ax(y)\|^2
	{}+{}
		(\conj g)^{\gamma}(y+\gamma Ax(y)).
\end{equation}

The AME enjoys several favorable properties, some of which we now summarize.
For any \(\gamma>0\), \(\psi_\gamma\) is (strictly) continuous over \(\R^m\), whereas if \(\gamma\) is small enough then the problem of minimizing \(\psi_\gamma\) is equivalent to solving \eqref{eq:DualProblem}.
These properties are listed in the next result.

\begin{thm}\label{thm:AMEbounds}
For any \(\gamma>0\), \(\psi_\gamma\) is a strictly continuous function on \(\R^m\) satisfying
\begin{enumerate}
	\item\label{it:AMEleq}
		\(
			\psi_\gamma(y)
		{}\leq{}
			\psi(y) + \tfrac\gamma 2\|Ax(y) - z_\gamma(y)\|^2
		\),
	\item\label{it:AMEgeq}
		\(
			\psi_\gamma(y)
		{}\geq{}
			\psi(T_\gamma(y))
			{}+{}
			\tfrac{\gamma}{2}\left(1-\tfrac{\gamma\|A\|^2}{\mu_f}\right)
			\|Ax(y) - z_\gamma(y)\|^2
		\),
\end{enumerate}
for any \(y\in\R^m\).
In particular, if \(\gamma<\nicefrac{\mu_f}{\|A\|^2}\), then the following also holds
\begin{enumerate}[resume]
	\item\label{it:EquivFBE}
		\(\inf\psi_\gamma=\inf\psi\)
		and
		\(\argmin\psi_\gamma=\argmin\psi\).
\end{enumerate}
	\begin{proof}
		Strict continuity of \(\psi_\gamma\) follows immediately by the expression \eqref{eq:EquivFBE}.
		\begin{proofitemize}
		\item\ref{it:AMEleq}:
			follows by \cref{lem:TwoPointsIneq} using \(w = y\).
		\item\ref{it:AMEgeq}:
			due to strong convexity of \(f\), \(\conj{f}\) has \(1/\mu_f\)-Lipschitz gradient, and consequently
			\ifieee
				\begin{align}
					\nonumber
						\conj{f}(-\trans{A}T_\gamma(y))
					&{}\leq{}
						\conj{f}(-\trans{A} y)-\innprod{ Ax(y)}{T_\gamma(y)-y}
					\\
					\nonumber
					&\phantom{{}\leq{}}
						+\tfrac{1}{2\mu_f}\|\trans{A}(T_\gamma(y)-y)\|^2
					\\
					\nonumber
					&{}={}
						\conj{f}(-\trans{A} y)
						-\gamma\innprod{ Ax(y)}{Ax(y)-z_\gamma(y)}
					\\
					&\phantom{{}={}}
						+\tfrac{\gamma^2}{2\mu_f}\|\trans{A}(Ax(y)-z_\gamma(y))\|^2. \label{eq:aux2}
				\end{align}
			\else
				\begin{align}
					\nonumber
						\conj{f}(-\trans{A}T_\gamma(y))
					&{}\leq{}
						\conj{f}(-\trans{A} y)-\innprod{ Ax(y)}{T_\gamma(y)-y}
						+ \tfrac{1}{2\mu_f}\|\trans{A}(T_\gamma(y)-y)\|^2
					\\
					\nonumber
					&{}={}
						\conj{f}(-\trans{A} y)
						-\gamma\innprod{ Ax(y)}{Ax(y)-z_\gamma(y)}
					\\
					&\phantom{{}={}}
						+\tfrac{\gamma^2}{2\mu_f}\|\trans{A}(Ax(y)-z_\gamma(y))\|^2. \label{eq:aux2}
				\end{align}
			\fi
			Combining \eqref{eq:AMEvsFBE} with \eqref{eq:aux2}:
			\ifieee
				\begin{align}
					\nonumber
						\psi_\gamma(y)
					&{}\geq{}
						\psi(T_\gamma(y)) - \tfrac{\gamma^2}{2\mu_f}\|\trans{A}(Ax(y)-z_\gamma(y))\|^2
					\\
					\nonumber
					&\phantom{{}\leq{}}
						+ \tfrac{\gamma}{2}\|Ax(y)-z_\gamma(y)\|^2
					\\
					&{}\geq{}
						\psi(T_\gamma(y)) + \tfrac{\gamma}{2}\left(1-\tfrac{\gamma\|A\|^2}{\mu_f}\right)\|Ax(y)-z_\gamma(y)\|^2. \nonumber 
				\end{align}
			\else
				\begin{align}
					\nonumber
						\psi_\gamma(y)
					&{}\geq{}
						\psi(T_\gamma(y)) - \tfrac{\gamma^2}{2\mu_f}\|\trans{A}(Ax(y)-z_\gamma(y))\|^2
						+ \tfrac{\gamma}{2}\|Ax(y)-z_\gamma(y)\|^2
					\\
					&{}\geq{}
						\psi(T_\gamma(y)) + \tfrac{\gamma}{2}\left(1-\tfrac{\gamma\|A\|^2}{\mu_f}\right)\|Ax(y)-z_\gamma(y)\|^2. \nonumber 
				\end{align}
			\fi
		\item\ref{it:EquivFBE}:
			easily follows combining \ref{it:AMEleq} and \ref{it:AMEgeq} with \(y = y_\star \in Y_\star\), in light of the dual optimality condition  \eqref{eq:criticalOptimal}.
		\qedhere
		\end{proofitemize}
	\end{proof}
\end{thm}

		\subsection{Analogy with the dual Moreau envelope}
			\Cref{thm:AMEvsFBE} highlights a clear connection between the augmented Lagrangian, the forward-backward envelope and the alternating minimization algorithm.
This closely resembles the one, first noticed by Rockafellar \cite{rockafellar1973dual,rockafellar1976augmented}, relating the augmented Lagrangian, the Moreau envelope and the \DEF{method of multipliers} (also known as \DEF{augmented Lagrangian method}) by Hestenes and Powell~\cite{hestenes1969multiplier,powell1969method}.
Consider the general linear equality constrained convex problem
\begin{equation}\label{eq:ECCP}\begin{aligned}
	\minimize_{z\in\R^k}\ & g(z) \\
	\stt\ & Bz = b,
\end{aligned}\end{equation}
where \(\func{g}{\R^m}{\Rinf}\) is proper, closed, convex, \(B\in\R^{m\times k}\) and \(b\in\R^m\). When applied to the dual of \eqref{eq:ECCP}, namely
\[
	\minimize_{y\in\R^m}\ \omega(y) = \conj g(-\trans{B} y) + \innprod by,
\]
the proximal minimization algorithm \cite[§5.2]{bertsekas2015convex} is equivalent to the following augmented Lagrangian method
\ifieee
	\begin{align*}\label{eq:MultipliersMethod}
		z^k
	{}={} &
		\argmin_{z\in\R^n} \set{g(z)+ \innprod{y^k}{Bz-b} + \tfrac{\gamma}{2}\|Bz-b\|^2}
	\\
		y^{k+1}
	{}={} &
		y^k+\gamma(Bz^k-b).
	\end{align*}
\else
	\begin{align*}\label{eq:MultipliersMethod}
		z^k
	{}={} &
		{\argmin}_{z\in\R^n} \set{g(z)+ \innprod{y^k}{Bz-b} + \tfrac{\gamma}{2}\|Bz-b\|^2},
	\\
		y^{k+1}
	{}={} &
		y^k+\gamma(Bz^k-b).
	\end{align*}
\fi
If
\(
	\range(\trans{B})\cap\relint(\dom \conj g)
{}\neq{}
	\emptyset
\)
one can show, with a similar proof to that of \Cref{thm:AMEvsFBE}, that
the Moreau envelope of $\omega$ satisfies
\begin{align*}
	\omega^\gamma(y^k)
{}={} &
	-g(z^k) - \innprod{y^k}{Bz^k-b} - \tfrac\gamma2\|Bz^k-b\|^2
\\
{}={} &
	-\Ell_\gamma(z^k, y^k).
\end{align*}
Therefore the forward-backward and Moreau envelope functions have the same nice interpretation in terms of augmented Lagrangian, when they are applied to the dual of equality
constrained convex problems: in a sense, \Cref{thm:AMEvsFBE} extends and generalizes the classical result on the dual Moreau envelope, by allowing for an additional variable \(x\) and a strongly convex term \(f\) in the problem.

	\section{Convergence}
		\label{sec:Convergence}
		We now turn our attention to the global convergence properties of \Cref{alg:NAMA}.
In light of \Cref{rem:GeneralizesAMA}, the results in this section directly apply to AMA, which is a special case of NAMA.

\begin{rem}[Termination of line-search]\label{rem:LineSearch}
The line-search step \ref{step:AMLS1_linesearch} is well defined regardless of the choice of $d^k$: at any iteration \(k\), condition \eqref{eq:AMLS_ascent} holds for $i_k$ sufficiently large.
To see this, suppose that \(\|Ax^k - z^k\| > 0\) (otherwise \((x^k, y^k)\) is a primal-dual solution). Then, since \(\gamma < \nicefrac{\mu_f}{\|A\|^2}\), \Cref{thm:AMEbounds} implies that
\begin{equation}\label{eq:FeasibleLS}
	\psi_{\gamma}(T_\gamma(y^k)) < \psi_{\gamma}(y^k).
\end{equation}
Since \(\lspoint y^k \to T_\gamma(y^k)\) as \(\tau_k\to 0\) and \(\psi_{\gamma}\) is continuous, then necessarily \(\psi_{\gamma}(\lspoint y^k) \leq \psi_{\gamma}(y^k)\) for \(\tau_k\) sufficiently small.
\end{rem}

\begin{rem}[Bounded iteration complexity]
	\label{rem:IterationComplexity}
In the best case where \(\tau_k=1\) is accepted in step \ref{step:AMLS1_linesearch}, exactly two alternating minimizations are performed at iteration \(k\).
In practice, one can also impose a lower bound \(\tau_{\min} > 0\) for \(\tau_k\): when \(\tau_k < \tau_{\min}\) then the ordinary AMA update \(y^{k+1} = y^k + \gamma(Ax^k-z^k)\) is executed and the algorithm proceeds to the next iteration.
This strategy results in a bounded iteration complexity for \algacronym, and does not affect the convergence results of this and later sections.
\end{rem}


\Cref{thm:AMEbounds} ensures that the following chain of inequalities, which will be fundamental for convergence results, holds in \Cref{alg:NAMA}:
\begin{subequations}\label{subeq:chain}
\begin{align}
	\psi(y^{k+1})
{}\leq{}&
	\psi_\gamma(\lspoint y^k)
\label{eq:chain1} \\
{}\leq{}&
	\psi_\gamma(y^k)
\label{eq:chain2}
\\
\label{eq:chain3}
{}\leq{}&
\psi(y^k)
	- \tfrac{\gamma}{2}\|Ax^k - z^k\|^2.
\end{align}
\end{subequations}
In particular, \Cref{alg:NAMA} is a descent method for \(\psi\).


We now prove that the iterates of \eqref{alg:NAMA} converge to the dual optimal cost and to the primal solution. Moreover, global convergence rates are provided.

\begin{thm}[Global convergence]\label{thm:GlobalConvergence}
In \Cref{alg:NAMA}:
\begin{enumthm}
	\item\label{it:subseqConv}
		\(x^k\to x_\star\),
		\(z^k\to Ax_\star\),
		and all cluster points of \(\seq{y^k}\)
		are dual optimal, \ie, they belong to \(Y_\star\);
	\item\label{it:sublinConv}
		if \(0 \in \interior(\dom g - A\dom f)\)
		then \(\psi(y^k)\searrow\inf\psi\) with global rate \(O(1/k)\), and \(x^k\to x_\star\) with global rate \(O(1/\sqrt{k})\);
	\item\label{it:linConv}
		if \(f\) and \(g\) are piecewise linear-quadratic
		then \(\psi(y^k)\searrow\inf\psi\) with global Q-linear rate, and \(x^k\to x_\star\) with global R-linear rate.
\end{enumthm}
			\begin{proof}
			\begin{proofitemize}
			\item\ref{it:subseqConv}:
				by \eqref{eq:chain3}, for all \(i \geq 0\) we have
				\[\psi(y^{i+1}) \leq \psi(y^{i}) - \tfrac{\gamma}{2}\|Ax^i -  z^i\|^2.\]
				By summing up the inequality for \(i=1,\ldots, k\) we obtain
				\[
					\inf \psi \leq \psi(y^{k+1}) \leq \psi(y^1) - \frac\gamma 2\sum_{i=1}^k \|Ax^i -  z^i\|^2
				\]
				(the sum starts from \(i=1\) since \(y^0\) may be dual infeasible).
				In particular (cf. \eqref{eq:TR}) \(R_\gamma(y^k) = z^k - Ax^k \to 0\), and since \(R_\gamma\) is continuous, necessarily all cluster points of \(\seq{y^k}\) are optimal.
				Moreover, it follows from \cref{lem:dual2primal} that the sequence \(\seq{x^k}\) is bounded.
				Let \(K\subseteq\N\) and \(\bar x\) be such that \(\seq{x^k}[k\in K]\to\bar x\); then, since \(Ax^k-z^k\to0\) we also have that \(\seq{z^k}[k\in K]\to A\bar x\).
				By multiplying \eqref{eq:OptCondy} on the left by \(\trans A\) and summing \eqref{eq:OptCondx} we obtain
				\(
					\gamma\trans A
					(Ax_k-z_k)
				{}\in{}
					\partial f(x_k)
					{}+{}
					\trans A\partial g(z_k)
				\).
				By letting \(K\ni k\to\infty\), from outer semicontinuity of the subdifferential we obtain that
				\[
					0
				{}\in{}
					\partial f(\bar x)
					{}+{}
					\trans A\partial g(A\bar x)
				{}\subseteq{}
					\partial(f+g\circ A)(\bar x)
				\]
				where the last inclusion follows from \cite[Thm.s 23.8 and 23.9]{rockafellar1997convex}.
				Thus, \(\bar x\) is optimal, and being \(x_\star\) the unique primal optimal (due to strong convexity), necessarily \(\bar x=x_\star\).
				From the arbitrarity of the cluster point we conclude that \(x^k\to x_\star\) and \(z^k\to Ax_\star\).
			\item\ref{it:sublinConv}:
				the assumed condition is equivalent to \(Y_\star\) being nonempty and compact, see \cite[Thm. 5.2.1]{auslender2003asymptotic}, which implies that \(\psi\) has bounded level sets \cite[Prop. 3.23]{rockafellar2011variational}. The proof proceeds similarly to that of~\cite[Thm. 4]{nesterov2013gradient}.
				Let \(D > 0\) be such that \(\dist(y,Y_\star)<D\) for all points \(y\in \set{y\in\R^m}[\psi(y) \leq \psi(y^0)]\).
				From \cite[Prop. 2.5]{stella2017forward} we know that \(\psi_\gamma \leq \psi^\gamma\) (the Moreau envelope of \(\psi\)).
				Therefore,
				\[
					\mathtight
					\psi(y^{k+1})
				{}\overrel[\leq]{\eqref{eq:chain2}}{}
					\psi_\gamma(y^k)
				{}\leq{}
					\psi^\gamma(y^k)
				{}={}
					\min_{\mathclap{w\in\R^m}}\set{\psi(w)+\tfrac{1}{2\gamma}\|w-y^k\|^2}
				\]
				and in particular, for $y_\star\in\argmin\psi$,
				{\mathtight\begin{align*}
					\psi(y^{k+1})
				{}\leq{} &
					\min_{
							\alpha\in[0,1]
					}{
						\set{
							\psi(\alpha y_\star+(1-\alpha)y^k)
							{}+{}
							\tfrac{\alpha^2}{2\gamma}\|y^k-y_\star\|^2
						}
					}
				\\
				{}\leq{} &
					\min_{
							\alpha\in[0,1]
					}{
						\set{
							\psi(y^k)-\alpha(\psi(y^k)-\inf\psi)+\tfrac{D^2}{2\gamma}\alpha^2
						}
					}
				\end{align*}}%
				where in last inequality we used convexity of $\psi$.
				In case $\psi(y^0)-\inf\psi\geq D^2/\gamma$, then the optimal solution of the latter problem for $k=0$ is $\alpha=1$, and
				\(
					\psi(y^1) - \inf\psi
				{}\leq{}
					\nicefrac{D^2}{2\gamma}
				\).
				Otherwise, the optimal solution is
				\[
					\alpha
				{}={}
					\tfrac{\gamma}{D^2}
					(\psi(y^k)-\inf\psi)
				{}\leq{}
					\tfrac{\gamma}{D^2}
					(\psi(y^0)-\inf\psi)
				{}\leq{}
					1
				\]
				and we obtain
				\[
					\psi(y^{k+1})
				{}\leq{}
					\psi(y^k)
					{}-{}
					\tfrac{\gamma}{2D^2}
					(\psi(y^k)-\inf\psi)^2.
				\]
				By letting $\lambda_k=\frac{1}{\psi(y^k)-\inf\psi}$ the last inequality becomes
				\[
					\lambda_{k+1}^{-1}\leq\lambda_k^{-1}-\tfrac{\gamma}{2D^2}\lambda_{k+1}^{-2}.
				\]
				By multiplying both sides by $\lambda_{k}\lambda_{k+1}$ and rearranging,
				\begin{align*}
					\lambda_{k+1}\geq\lambda_k+\tfrac{\gamma}{2D^2}\tfrac{\lambda_{k+1}}{\lambda_k}\geq \lambda_k+\tfrac{\gamma}{2D^2},
				\end{align*}
				where the latter inequality follows from the fact that the sequence $\seq{\psi(y^k)}$ is nonincreasing, as shown in \eqref{subeq:chain}.
				By telescoping the inequality we obtain
				\[
					\lambda_k
				{}\geq{}
					\lambda_0
					{}+{}
					k\tfrac{\gamma}{2D^2}
				{}\geq{}
					k\tfrac{\gamma}{2D^2},
				\]
				and therefore
				\(
					\psi(y^k) - \inf\psi
				{}\leq{}
					\nicefrac{2D^2}{k\gamma}
				\). This, together with \cref{lem:dual2primal}, proves \ref{it:sublinConv}.
			\item\ref{it:linConv}:
				since the primal optimum is finite (see \cref{rem:Basic}), if \(f\) and \(g\) are piecewise linear-quadratic then \(Y_\star\) is nonempty, see \cite[Thm. 11.42, Ex. 11.43]{rockafellar2011variational}.
				Using \eqref{subeq:chain} we have that
				\begin{equation}\label{eq:LinearIneq2}
					\psi(y^k) - \psi(y^{k+1}) \geq \tfrac{\gamma}{2}\|Ax^k-z^k\|^2.
				\end{equation}
				Furthermore, using \cref{lem:TwoPointsIneq} with \(w = y_\star^k = \Pi_{Y_\star}y^k\) and \(y = y^k\), we obtain
				\begin{align*}
					\psi(y^{k+1}) - \inf\psi &{}\leq{} \psi_\gamma(y^k) - \inf\psi \\
					&{}\leq{} \innprod{Ax^k-z^k}{y_\star^k - y^k} - \tfrac{\gamma}{2}\|Ax^k-z^k\|^2,
				\end{align*}
				where first inequality is due to \eqref{eq:chain3}. This implies
				\[
					\psi(y^{k+1}) - \inf\psi \leq \|Ax^k-z^k\|^2\left(\tfrac{\dist(y^k, Y_\star)}{\|Ax^k-z^k\|}-\tfrac{\gamma}{2}\right)
				\]
				which, by using \eqref{eq:LinearIneq2}, yields
				\begin{equation}\label{eq:LinearIneq3}
					\psi(y^{k+1}) - \inf\psi \leq \left(1 - \tfrac{\gamma}{2}\tfrac{\|Ax^k-z^k\|}{\dist(y^k, Y_\star)}\right)(\psi(y^{k}) - \inf\psi).
				\end{equation}
				It follows from \cite[Thm. 11.14]{rockafellar2011variational} that \(\conj{f}\) and \(\conj{g}\) are convex piecewise linear-quadratic in this case, and so is \(\psi\). Therefore by \cite[Thm. 2.7]{li1995error} \(\psi\)
				enjoys the following quadratic growth condition: for any \(\nu>0\) there is \(\alpha>0\) such that
				\[ \tfrac{\alpha}{2}\dist^2(y, Y_\star) \leq \psi(y^k)-\inf\psi\quad \forall y : \psi(y)-\inf\psi \leq \nu, \]
				which by \cite[Cor. 3.6]{drusvyatskiy2016error} is equivalent to the following error bound condition for some \(\beta > 0\)
				\begin{equation}\label{eq:ErrorBound}
					\dist(y, Y_\star) \leq \beta\|Ax(y)-z_\gamma(y)\|
				\end{equation}
				holding for all \(y\) such that \(\psi(y)-\inf\psi \leq \nu\).
				By using \eqref{eq:ErrorBound} in \eqref{eq:LinearIneq3} we obtain global Q-linear convergence of \(\seq{\psi(y^k)}\), and from \cref{lem:dual2primal} global R-linear convergence of \(\seq{x^k}\) also follows.
				\qedhere
		\end{proofitemize}
	\end{proof}
\end{thm}

In general
we can prove local linear convergence of \Cref{alg:NAMA} provided that \(\partial f\) and \(\partial g\) are \emph{calm}, according to the following definition (see \cite[Sec. 3H, Ex. 3H.4]{dontchev2009implicit}).

\begin{defin}[Calmness of a mapping]\label{def:Calmness}
	A multi-valued mapping \(F:\R^m\rightrightarrows\R^n\) is said to be \emph{calm} at \(\bar y\in\R^m\) for \(\bar x\in F(\bar y)\) if there is a neighborhood \(U\) of \(\bar x\) such that
	\[F(y)\cap U \subseteq F(\bar y) + O(\|y-\bar y\|), \quad \forall y\in\R^m.\]
	We simply say that \(F\) is calm at \(\bar y\in\R^m\) (with no mention of \(\bar x\)) if it is calm at \(\bar y\in\R^m\) for all \(\bar x\in F(\bar y)\).
\end{defin}

Calmness is a very common property of the subdifferential mapping.
The subdifferential of all piecewise linear-quadratic functions is calm everywhere, as follows from \cite[Prop. 3H.1]{dontchev2009implicit}. Other examples include the nuclear and spectral norms \cite{schoepfer2016linear}.
Smooth functions, \ie with Lipschitz gradient, clearly have calm subdifferential: this includes Moreau envelopes of closed, convex functions, such as the Huber loss for robust estimation, and commonly used loss functions such as the squared Euclidean norm and the logistic loss.

Calmness is equivalent to metric subregularity of the inverse mapping \cite[Thm. 3H.3]{dontchev2009implicit}: from \cite[Prop. 6, Prop. 8]{zhou2015unified} we then deduce that the indicator functions of \(\ell_1\), \(\ell_{\infty}\) and Euclidean norm balls all have calm subdifferentials.

The following result holds. Its proof is analogous to the one of \cite[Thm. 4.2]{drusvyatskiy2016error}, although our assumption of calmness is equivalent to metric subregularity of \(\partial \conj f\) and \(\partial \conj g\), which is implied by the firm convexity assumed in \cite{drusvyatskiy2016error}.

\begin{thm}[Local linear convergence]\label{thm:LinearConvergence}
Suppose that the following hold for \eqref{eq:Problem}:
\begin{enumerate}
	\item\label{it:ThmStrongFeas} \(0 \in \interior(\dom g - A\dom f)\) (nonempty, compact \(Y_\star\));
	\item\label{it:ThmStrictComp} \(0\in\relint\partial(f+g\circ A)(x_\star)\) (strict complementarity).
\end{enumerate}
Suppose also that \(\partial f\) is calm at \(x_\star\) and \(\partial g\) is calm at \(Ax_\star\).
Then in \Cref{alg:NAMA} eventually \(\psi(y^k) \to \inf\psi\) with \(Q\)-linear rate and \(x^k\to x_\star\) with \(R\)-linear rate.
\begin{proof}
	As discussed in the proof of \cref{it:linConv}, it suffices to show that an error bound of the form \eqref{eq:ErrorBound} holds for some \(\beta,\nu > 0\).

	The assumed calmness properties of \(\partial f\) and \(\partial g\) are equivalent to \emph{metric subregularity} of \(\partial \conj f\) at \(-\trans{A}y_\star\) for \(x_\star\), and of \(\partial \conj g\) at \(y_\star\) for \(Ax_\star\), see \cite[Thm. 3H.3]{dontchev2009implicit}, for all \(y_\star\in Y_\star\).
	This can be seen, using \cite[Thm. 3.3]{aragon2008characterization}, to be equivalent to the following: there exist \(c_{y_\star}>0\) and a neighborhood \(U_{y_\star}\) of \(y_\star\) such that for all \(y\in U_{y_\star}\)
	\ifieee
		\begin{align*}
			\conj f(-\trans{A}y)
		{}\geq{}&
			\phantom{{}+{}}
			\conj f(-\trans{A}y_\star) + \innprod{x_\star}{\trans{A}(y_\star-y)}
		\\
		&
			{}+\tfrac{c_{y_\star}}{2}\dist^2(-\trans{A}y, (\nabla\conj f)^{-1}(x_\star)),
		\\
			\conj g(y)
		{}\geq{}&
			\phantom{{}+{}}
			\conj g(y_\star) + \innprod{Ax_\star}{y - y_\star}
		\\
		&
			{}+ \tfrac{c_{y_\star}}{2}\dist^2(y, (\partial \conj g)^{-1}(Ax_\star)).
		\end{align*}
	\else
		\begin{align*}
			\conj f(-\trans{A}y)
		{}\geq{}&
			\conj f(-\trans{A}y_\star) + \innprod{x_\star}{\trans{A}(y_\star-y)}
			+\tfrac{c_{y_\star}}{2}\dist^2(-\trans{A}y, (\nabla\conj f)^{-1}(x_\star)),
		\\
			\conj g(y)
		{}\geq{}&
			\conj g(y_\star) + \innprod{Ax_\star}{y - y_\star}
			+ \tfrac{c_{y_\star}}{2}\dist^2(y, (\partial \conj g)^{-1}(Ax_\star)).
		\end{align*}
	\fi
	Since \(Y_\star \subset \bigcup_{y_\star\in Y_\star}U_{y_\star}\) and \(Y_\star\) is nonempty and compact (due to \ref{it:ThmStrongFeas}, see \cite[Thm. 5.2.1]{auslender2003asymptotic}), we may select a \emph{finite} subset \(W\subset Y_\star\) such that \(Y_\star\subset U_{Y_\star} = \bigcup_{y_\star\in W}U_{y_\star}\).
	Summing the above inequalities for all \(y_\star\in W\), and denoting \(c = \min\set{c_{y_\star}}[y_\star\in W] > 0\), we obtain
	\ifieee
		\begin{multline}
			\psi(y)
		{}\geq{}
			\inf\psi
		\\\label{eq:PreQGC}
			{}+{}
			\tfrac{c}{2}\left[\dist^2(-\trans{A}y, \partial f(x_\star)) + \dist^2(y, \partial g(Ax_\star))\right]
		\end{multline}
	\else
		\begin{equation}\label{eq:PreQGC}
			\psi(y)
		{}\geq{}
			\inf\psi + \tfrac{c}{2}\left[\dist^2(-\trans{A}y, \partial f(x_\star)) + \dist^2(y, \partial g(Ax_\star))\right],\quad \forall y\in U_{Y_\star},
		\end{equation}
	\fi
	for all \(y\in U_{Y_\star}\), where we have also used \((\nabla\conj f)^{-1} = \partial f\) and \((\partial \conj g)^{-1} = \partial g\).
	Note that \ref{it:ThmStrongFeas} implies strict feasibility, therefore from \cref{lem:linearReg}, and the fact that for any \(a,b\in\R\), \(a^2+b^2 \geq 2ab\), we obtain that \eqref{eq:PreQGC} implies
	\[\psi(y) \geq \inf\psi + \tfrac{\kappa}{2}\dist^2(y, Y_\star),\quad \forall y\in U_{Y_\star}, \]
	for some \(\kappa > 0\), \ie \(\psi\) satisfies the quadratic growth condition, which by \cite[Cor. 3.6]{drusvyatskiy2016error} is equivalent to the error bound condition \eqref{eq:ErrorBound}.
	This completes the proof.
\end{proof}
\end{thm}

\begin{rem}[Backtracking on \(\gamma\)]
In practice, no prior knowledge of the global Lipschitz constant \(\nicefrac{\|A\|^2}{\mu_f}\) is required for \Cref{alg:NAMA}:
instead of a fixed parameter \(\gamma\), one can adaptively determine a \emph{sequence} \(\seq{\gamma_k}\) essentially ensuring that inequalities \eqref{eq:FeasibleLS} (which guarantees termination of the line-search step \ref{step:AMLS1_linesearch}) and \eqref{eq:chain1} (which guarantees descent) hold at every iteration.
This is done as follows. Select \(\alpha\in(0,1)\) and initialize \(\gamma_0 > 0\).
At iteration \(k\), let \(\bar y^k = y^k + \gamma_k(Ax^k - z^k)\) and \(\bar x^k = x(\bar y^k)\), and if
\[
	f(x^k) > f(\bar x^k) - \innprod{\trans{A}\bar y^k}{x^k - \bar x^k} + \tfrac{\alpha\gamma}{2}\|Ax^k - z^k\|^2,
\]
then \(\gamma_k \gets \gamma_k/2\) and restart the iteration. Similarly if
\[
	f(\tilde x^k) > f(x^{k+1}) - \innprod{\trans{A}y^{k+1}}{\tilde x^k - x^{k+1}} + \tfrac{\alpha\gamma}{2}\|A\tilde x^k - \tilde z^k\|^2.
\]
As soon as \(\gamma_k \leq \nicefrac{\alpha\mu_f}{\|A\|^2}\), the two inequalities above will never hold. As a consequence, \(\gamma_k\) will be decreased only a finite number of times and will be constant starting from some iteration \(\bar k\).
The inequalities above are obtained by imposing the usual quadratic upper bound on \(\conj{f}\circ(-\trans{A})\), due to smoothness, and applying the conjugate subgradient theorem \eqref{eq:ConjSubgr} in light of \eqref{eq:OptCondx}.
This procedure of adaptively adjusting \(\gamma_k\) is analogous to what is done in practice in (fast) AMA, see \cite[Rem. 3.4]{beck2014fast} and \cite[§3, §4]{beck2009fast},
and does not affect the validity of \cref{thm:GlobalConvergence,thm:LinearConvergence}.
\end{rem}

	\section{First- and second-order properties}
		\label{sec:DiffFBE}
\Cref{alg:NAMA} is a line-search method for the unconstrained minimization of \(\psi_\gamma\) which, by \Cref{it:EquivFBE}, is equivalent to solving \eqref{eq:DualProblem}.
To enable fast convergence of the iterates, we can apply ideas from smooth unconstrained optimization in selecting the sequence \(\seq{d^k}\) of directions. To this end, differentiability of \(\psi_\gamma\) around dual solutions \(y_\star\) is a desirable property.
We will now see that this is implied by generalized second-order properties of \(f\) around \(x_\star\), which are introduced in the following assumption.
Analogous assumptions on \(g\) further ensure that \(\psi_\gamma\) is (strictly) twice differentiable at \(y_\star\).
The interested reader is referred to \cite{rockafellar2011variational} for an extensive discussion on (second-order) epi-differentiability.

\begin{ass}
	\label{ass:fg2}
The following hold with respect to a primal-dual solution \((x_\star,y_\star)\) to \eqref{eq:Problem}-\eqref{eq:DualProblem}:
\begin{enumerate}
\item\label{ass:f2}
	\(f\) is strictly twice epi-differentiable at all \(x\in\dom f\) close enough to \(x_\star\),
	and in particular the second-order epi-derivative at \(x_\star\) for \(-\trans{A} y_\star\)
	\ifieee
		is, for \(w\in\R^n\),
		\begin{equation}\label{eq:fGenQuad}
			\twiceepi[-\trans{A} y_\star]{f}{x_\star}[w]
		{}={}
			\innprod{H_fw}{w} + \delta_{S_f}(w),
		\end{equation}
	\else
		is
		\begin{equation}\label{eq:fGenQuad}
			\twiceepi[-\trans{A} y_\star]{f}{x_\star}[w]
		{}={}
			\innprod{H_fw}{w} + \delta_{S_f}(w),
		\quad\forall w\in \R^m,
		\end{equation}
	\fi
	where \(S_f\) is a linear subspace of \(\R^n\) and \(H_f\in\R^{n\times n}\);
\item\label{ass:g2}
	\(g\) is (strictly) twice epi-differentiable at \(A x_\star\) for \(y_\star\),
	\ifieee
		with
		\begin{equation}\label{eq:GenQuad}
			\twiceepi[y_\star]{g}{Ax_\star}[w]
		{}={}
			\innprod{H_gw}{w} + \delta_{S_g}(w),
		\end{equation}
		for all \(w\in\R^m\),
	\else
		with
		\begin{equation}\label{eq:GenQuad}
			\twiceepi[y_\star]{g}{Ax_\star}[w]
		{}={}
			\innprod{H_gw}{w} + \delta_{S_g}(w),
			\quad\forall w\in \R^m,
		\end{equation}
	\fi
	where \(S_g\) is a linear subspace of \(\R^m\) and \(H_g\in\R^{m\times m}\).
\end{enumerate}
When the stronger condition in parenthesis holds we will say that the assumptions are \emph{strictly} satisfied.

Without loss of generality, we consider \(H_f\) and \(H_g\) symmetric and positive semidefinite, satisfying \(\range(H_f)=S_f\), \(\nullspace(H_f)=S_f^\bot\), \(\range(H_g)\subseteq S_g\) and \(\nullspace(H_g)\supseteq S_g^\bot\).
\end{ass}
The requirements on \(H_f\) and \(H_g\) can indeed be made without loss of generality: matrix
\(
	H_f'
{}={}
	\frac 12\proj_{S_f}(H_f+\trans H_f)\proj_{S_f}
\)
has the desired properties and satisfies \eqref{eq:fGenQuad} provided \(H_f\) does, and similarly for \(H_g\).
In particular, it holds that
\begin{equation}\label{eq:H}
	H_f
{}={}
	\proj_{S_f}H_f\proj_{S_f}
\quad\text{and}\quad
	H_g
{}={}
	\proj_{S_g}H_g\proj_{S_g}.
\end{equation}
\begin{thm}[Differentiability of \(\psi_\gamma\)]\label{thm:FirstOrder}
Suppose that \Cref{ass:f2} holds for a primal-dual solution \((x_\star,y_\star)\).
Then \(\psi_\gamma\) is of class \(\cont^1\) around \(y_\star\), with
\[
	\nabla\psi_\gamma(y)
{}={}
	Q_\gamma(y)R_\gamma(y)
\]
where \(Q_\gamma(y) = I-\gamma A\nabla^2\conj f(-\trans{A} y)\trans{A}\).
\begin{proof}
From \cref{lem:psi1TwiceDiff} it follows that \(\tf=\conj f\circ(-\trans A)\) is of class \(\cont^2\) around \(y_\star\).
The claim now easily follows from the chain rule of differentiation applied to \eqref{eq:EquivFBE}, by using \eqref{eq:MoreauGradient}.
\end{proof}
\end{thm}

		Twice differentiability of \(\psi_\gamma\) at a dual solution \(y_\star\) is very important: when Newton-type directions are used, this implies that eventually unit stepsize will be accepted and fast asymptotic convergence will take place. In other words, unlike standard nonsmooth merit functions for constrained optimization, \(\psi_\gamma\) does not prevent the acceptance of unit stepsize.

\begin{thm}[Twice differentiability of \(\psi_\gamma\)]\label{thm:TwiceDiff}
	Suppose that \Cref{ass:fg2} (strictly) holds with respect to a primal-dual solution \((x_\star,y_\star)\).
	Then,
	\begin{enumerate}
		\item\label{it:Rdiff}
			\(R_\gamma\) is (strictly) differentiable at \(y_\star\) with Jacobian
			\begin{equation}\label{eq:JacFPR}
				\jac{R_\gamma}{y_\star} = \gamma^{-1}\left[I-P_\gamma(y_\star)Q_\gamma(y_\star)\right];
			\end{equation}
			here, \(Q_\gamma\) is as in \Cref{thm:FirstOrder} and
			\ifieee
				\begin{align}
				\nonumber
					P_\gamma(y_\star)
				{}={} &
					J\prox_{\gamma\conj g}\left(y_\star+\gamma A\nabla\conj f(-\trans Ay_\star)\right)
				\\
				\label{eq:JacP}
				{}={} &
					\proj_{\bar S}
					\left(I+\gamma\geninv{H_g}\right)^{-1}
					\proj_{\bar S}
				\end{align}
			\else
				\begin{equation}\label{eq:JacP}
					P_\gamma(y_\star)
				{}={}
					J\prox_{\gamma\conj g}\left(y_\star+\gamma A\nabla\conj f(-\trans Ay_\star)\right)
				{}={}
					\proj_{\bar S}
					\left(I+\gamma\geninv{H_g}\right)^{-1}
					\proj_{\bar S}
				\end{equation}
			\fi
			with \(\bar S={S_g}^\bot + \range(H_g)\);
		\item\label{it:psigammaTwiceDiff}
			\(\psi_\gamma\) is (strictly) twice differentiable at \(y_\star\) with symmetric Hessian%
			\begin{equation}\label{eq:HessianAME}
				\nabla^2\psi_\gamma(y_\star)
			{}={}
				\gamma^{-1}
				Q_\gamma(y_\star)
				\bigl[
					I-P_\gamma(y_\star) Q_\gamma(y_\star)
				\bigr].
			\end{equation}
	\end{enumerate}
	\begin{proof}
		Let \(\tf=\conj f\circ(-\trans A)\) and \(L_\tf=\nicefrac{\mu_f}{\|A\|^2}\).
		We know from \cite[Thms. 3.8, 4.1]{poliquin1996generalized} and \cite[Thm. 13.21]{rockafellar2011variational} that \(\prox_{\gamma\conj g}\) is (strictly) differentiable at \(\Fw[\tf]{y_\star}\) if and only if \(g\) (strictly) satisfies \cref{ass:g2}; in fact, by \eqref{eq:FirstOrderNecessary} we know that \(Ax_\star=-\nabla\tf(y_\star)\).
		Moreover, due to \cref{lem:psi1TwiceDiff}, \(\tf\in \cont^2\) in a neighborhood of \(y_\star\) and in particular $\nabla\tf$ is strictly differentiable at \(y_\star\).
		The formula for \(\jac{R_\gamma}{y_\star}\) follows from \eqref{eq:MoreauGradient} and the chain rule of differentiation.

		We now prove the claimed expression for $P_\gamma(y_\star)$.
		We may invoke \cref{lem:TwiceEpiDuality} and apply \cite[Ex. 13.45]{rockafellar2011variational} to the \emph{tilted} function $g+\innprod{\nabla\tf(y_\star)}{{}\cdot{}}$
		which this tells us that
		for all \(d\in\R^m\)
		\ifieee
			\begin{align*}
			&
				P_\gamma(y_\star)d
			\\
			={}&
				\prox_{(\gamma/2)\twiceepi[Ax_\star]{\conj g}{y_\star}}(d)
			\\
			={}&
				\argmin_{d'\in\bar S}{
					\set{
						\tfrac12\innprod{d'}{\geninv{H_g}d'} + \tfrac{1}{2\gamma}\|d' - d\|^2
					}
				}
			\\
			={}&
				\proj_{\bar S}
				\argmin_{d'\in\R^n}{
					\set{
						\tfrac12\innprod{\proj_{\bar S} d'}{\!\!\geninv{H_g}\proj_{\bar S} d'} + \tfrac{1}{2\gamma}\|\proj_{\bar S} d' - d\|^2
					}
				}
			\\
			={}&
				\proj_{\bar S}
				 {\bigl(
					\proj_{\bar S}[I+\gamma\geninv{H_g}]\proj_{\bar S}
				\bigr)}^\dagger
				  \proj_{\bar S} d
			\end{align*}
		\else
			\begin{align*}
				P_\gamma(y_\star)d
			{}={}&
				\prox_{(\gamma/2)\twiceepi[Ax_\star]{\conj g}{y_\star}}(d)
			\\
			{}={}&
				\argmin_{d'\in\bar S}{
					\set{
						\tfrac12\innprod{d'}{\geninv{H_g}d'} + \tfrac{1}{2\gamma}\|d' - d\|^2
					}
				}
			\\
			{}={}&
				\proj_{\bar S}
				\argmin_{d'\in\R^n}{
					\set{
						\tfrac12\innprod{\proj_{\bar S} d'}{\geninv{H_g}\proj_{\bar S} d'} + \tfrac{1}{2\gamma}\|\proj_{\bar S} d' - d\|^2
					}
				}
			\\
			{}={}&
				\proj_{\bar S}
				{\bigl(
					\proj_{\bar S}[I+\gamma\geninv{H_g}]\proj_{\bar S}
				 \bigr)}^\dagger
				\proj_{\bar S} d
			\\
			{}={}&
				\proj_{\bar S}
				[I+\gamma\geninv{H_g}]^{-1}
				\proj_{\bar S} d
			\end{align*}
		\fi
		where $^\dagger$ indicates the pseudo-inverse.
		Observe now that, since \(\range\geninv{H_g}=\range H_g\subseteq\bar S\), we have
		\[
			\proj_{\bar S}[I+\gamma\geninv{H_g}]\proj_{\bar S}
		{}={}
			AB
		\quad\mbox{for}\quad
			A
		{}={}
			I+\gamma\geninv{H_g}~\text{and}~
			B
		{}={}
			\proj_{\bar S}.
		\]
		Moreover,
		\begin{align*}
			\range(\trans AAB)
		{}\subseteq{} &
			\range B,
		\\
			\range(\trans BBA)
		{}\subseteq{} &
			\R^n=\range(A),
		\end{align*}
		therefore we can apply \cite[Facts 6.4.12 (i)-(ii) and 6.1.6 (xxxii)]{bernstein2009matrix} to see that
		\(
			{\bigl(
			\proj_{\bar S}[I+\gamma\geninv{H_g}]\proj_{\bar S}
		\bigr)}^\dagger
		{}={}
			\proj_{\bar S}
			[I+\gamma\geninv{H_g}]^{-1}
		\),
		yielding \eqref{eq:JacP}.

		Since $R_\gamma(y_\star) = 0$, from \cite[Lem. 6.2]{stella2017forward} it follows that \(\nabla\psi_\gamma=Q_\gamma R_\gamma\) is (strictly) differentiable at \(y_\star\) provided that $Q_\gamma$ is (strictly) continuous at $y_\star$ and $R_\gamma$ is (strictly) differentiable at $y_\star$.
		A simple application of the chain rule of differentiation concludes the proof of \ref{it:psigammaTwiceDiff}.
		\end{proof}
\end{thm}

To better understand the requirements of \Cref{ass:fg2}, let us consider the following simple but significant example: when \(f\) is \(\cont^2\) and \(g\circ A\) models linear inequality constraints, \Cref{ass:fg2} is implied by strict complementarity.

\begin{es}[\(\cont^2\) functions subject to polyhedral constraints]\label{ex:QP} Consider problems of the form
\[
	\minimize_{x\in\R^n}\ f(x) + \delta_C(Ax), \label{eq:QP}
\]
where 
\(g = \delta_C\) is the indicator of \(C = \set{z\in\R^m}[z\leq b]\), \(b\in\R^m\), and \(f\in \cont^2\).
In this case \Cref{ass:f2} holds with \(H_f = \nabla^2 f(x_\star)\),
\(S_f = \R^n\) (therefore \(\proj_{S_f} = \id\)), see \cite[Ex. 13.8]{rockafellar2011variational}.
Regarding \Cref{ass:g2}, one can use \cite[Ex. 13.17]{rockafellar2011variational} to see that
\[
	\twiceepi[y_\star]{g}{Ax_\star}[w]
{}={}
	\delta_{K(Ax_\star,y_\star)}(w),
\]
where \(K\) is the critical cone. Denoting as \(T_C(y)\) the tangent cone of set \(C\) at \(y\in C\), and as \(J=\set{i}[(Ax_\star)_i = b_i]\) the set of active constraints at the solution \(x_\star\), the critical cone is given by
\begin{align*}
	K(Ax_\star,y_\star) &{}={} \set{w\in T_C(Ax_\star)}[\innprod{y_\star}{w} = 0]\\
	&{}={} \set{w}[\innprod{y_\star}{w} = 0, w_i \leq 0\;\forall i\in J].
\end{align*}
For \(K(Ax_\star,y_\star)\) to be a subspace, necessarily \((y_\star)_i > 0\) for all \(i\in J\), \ie, strict complementarity must hold at the primal-dual solution \((x_\star,y_\star)\).
In this case, \Cref{ass:g2} holds with \(H_g = 0\) and
\[
	S_g = K(Ax_\star,y_\star) = \set{w}[w_i = 0\;\forall i\in J].
\]
We may assume that \(J=\set{1,\ldots,k}\) without loss of generality, \ie, the first \(k\) constraints are the active ones, and let \(\bar J = \set{1,\ldots,m}\setminus J\). Note that \(\nabla^2\conj{f}(-\trans{A}y_\star) = \nabla^2 f(x_\star)^{-1}\) due to strong convexity of \(f\), see \cite[Ex. 11.9]{rockafellar2011variational}. By partitioning the inverse Hessian and constraint matrix as
\[
\nabla^2 f(x_\star)^{-1} = \begin{bmatrix}
H_{JJ} &\quad H_{J\bar J} \\
H_{\bar JJ} &\quad H_{\bar J\bar J} \\
\end{bmatrix},\quad
A = \begin{bmatrix}
A_J \\
A_{\bar J}
\end{bmatrix},\]
and using the notation of \Cref{it:Rdiff} we obtain
\[
	P_\gamma(y_\star) =
		\begin{bmatrix}I_k & 0 \\ 0 & 0 \end{bmatrix},
	\quad
	\jac{R_\gamma}{y_\star} =
		\begin{bmatrix}
			A_J H_{J\!J} \trans{A_J} & A_J H_{J\!\bar J} \trans{A_{\bar J}} \\
			0 & \tfrac1\gamma I_{m-k}
		\end{bmatrix},
\]
as it follows by elementary computations.
\end{es}

Finally, we can relate strong minimality of \(\psi\) and \(\psi_\gamma\) to nonsingularity of the Jacobian of \(R_\gamma\) and to the generalized second-order properties of \(f\) and \(g\) as follows.

\begin{thm}[Conditions for strong minimality]\label{thm:StrongMinim}
	If \Cref{ass:fg2} holds for a primal-dual solution \((x_\star,y_\star)\), then for all \(\gamma < \mu_f/\|A\|^2\) the following are equivalent:
	\begin{enumerateq}
		\item\label{thm:StrongMinimPsi}
			\(y_\star\) is a strong minimum for \(\psi\);%
			\footnote{We say that \(y_\star\) is a strong local minimum for \(h\) if for some \(\alpha>0\), \(\alpha\|y - y_\star\|^2 \leq h(y)-h(y_\star)\) for all \(y\) sufficiently close to \(y_\star\).}
		\item\label{thm:HessDef+}
			\(\nabla^2\psi_\gamma(y_\star)\) is nonsingular (in fact, positive definite);
		\item\label{thm:JRnonsing}
			\(\jac{R_\gamma}{y_\star}\) is nonsingular (in fact, similar to a symmetric and positive definite matrix);
		\item\label{thm:StrongMinimFBE}
			\(y_\star\) is a strong minimum for \(\psi_\gamma\).
	\end{enumerateq}
	\begin{proof}
		\begin{proofitemize}
		\renewcommand\myVar[3][\Leftrightarrow]{%
			\hspace*{0pt}%
			\mbox{%
				\textfillwidthof[c]{\ref*{thm:StrongMinimPsi}}{\ref{#2}}%
				~\(\fillwidthof[c]{\Leftrightarrow}{#1}\)~%
				\textfillwidthof[c]{\ref*{thm:StrongMinimPsi}}{\ref{#3}}:%
			}%
		}%
		\item\myVar{thm:HessDef+}{thm:JRnonsing}
			Let $P = P_\gamma(y_\star)$ and $Q = Q_\gamma(y_\star)$ for brevity.
			Notice first that, due to \cref{it:EquivFBE}, \(y_\star\) minimizes \(\psi_\gamma\) and therefore \(\nabla^2\psi_\gamma(y_\star)\succeq 0\).
			Moreover, since \(Q\) is symmetric and positive definite,
			\[
				\jac{R_\gamma}{y_\star}
			{}={}
				\gamma^{-1}(I-PQ)
			{}\sim{}
				Q^{-\nicefrac12}\nabla^2\psi_\gamma(y_\star)Q^{-\nicefrac12}
			\]
			the latter matrix being symmetric and positive semidefinite, where \(\sim\) denotes the similitude relation.
		\item\myVar{thm:HessDef+}{thm:StrongMinimFBE}
			trivial since $\nabla^2\psi_\gamma(y_\star)$ exists.
		\item\myVar{thm:StrongMinimFBE}{thm:StrongMinimPsi}
			the right implication is trivial since \(\psi_\gamma\leq\psi\) and \(\psi_\gamma(y_\star)=\psi(y_\star)\) as it follows from \cref{thm:AMEbounds}.
			Suppose now that there exist \(c,\varepsilon>0\) such that
			\(
				\psi(y)-\psi(y_\star)
			{}\geq{}
				\tfrac c2\|y-y_\star\|^2
			\)
			for all \(y\in\ball{y_\star}{\varepsilon}\).
			Since \(\conj g\) is convex, it follows that \(\prox_{\gamma\conj g}\) is \(1\)-Lipschitz continuous; combined with the fact that \(\nabla\conj f\) is \(\tfrac{1}{\mu_f}\)-Lipschitz continuous, we obtain that the alternating minimization operator \(T_\gamma\) is Lipschitz continuous with modulus \(\nicefrac{\|A\|^2}{\mu_f}\).
			Let \(\varepsilon'=\nicefrac{\mu_f}{\|A\|^2}\varepsilon\); since \(T_\gamma(y_\star)=y_\star\), for all \(y\in\ball{y_\star}{\varepsilon'}\) necessarily \(T_\gamma(y)\in\ball{y_\star}{\varepsilon}\).
			Therefore, letting
			\(
				c'
			{}={}
				\min\set{
					c,\,
					\gamma
					\bigl(1-\tfrac{\gamma\|A\|^2}{\mu_f}\bigr)
				}
			{}>{}
				0
			\),
			it follows from \cref{it:AMEgeq} that for all \(y\in\ball{y_\star}{\varepsilon'}\)
			\begin{align*}
				\psi_\gamma(y)
				{}-{}
				\psi_\star
			{}\geq{} &
				\psi(T_\gamma(y))
				{}-{}
				\psi_\star
				{}-{}
				\tfrac{\gamma}{2}
				\bigl(1-\tfrac{\gamma\|A\|^2}{\mu_f}\bigr)
				\|y-T_\gamma(y)\|^2
			\\
			{}\geq{} &
				\tfrac{c'}2
				\left(
					\|T_\gamma(y)-y_\star\|^2
					{}+{}
					\|y-T_\gamma(y)\|^2
				\right)
			\\
			{}\geq{} &
				\tfrac{c'}4
				\|y-y_\star\|^2.
			\end{align*}
			This shows that \(y_\star\) is a strong local minimum for \(\psi_\gamma\).
		\qedhere
		\end{proofitemize}
	\end{proof}
\end{thm}

In the context of \Cref{ex:QP}, notice that
\begin{equation*}
\jac{R_\gamma}{y_\star}\ \mbox{is nonsingular}\ \Leftrightarrow\
A_J H_{JJ} \trans{A_J}\ \mbox{is nonsingular}.
\end{equation*}
Since \(\nabla^2 f(x_\star) \succ 0\) by assumption, then \(H_{JJ} \succ 0\) and nonsingularity of the Jacobian is equivalent to \(A_J\) being full row rank,
\ie, linear independence of the active constraints at \(x_\star\) (the LICQ assumption).

	\section{Superlinear convergence}
		\label{sec:Superlinear}
		The following definition (cf. \cite[Eq. (7.5.2)]{facchinei2003finite}) gives the fundamental condition, on the sequence \(\seq{d^k}\) of directions, ensuring superlinear asymptotic convergence of \Cref{alg:NAMA}.

\begin{defin}[Superlinear directions]\label{def:SuperlinearDir}
	For \(\seq{y^k}\) converging to \(y_\star\), we say that \(\seq{d^k}\) is \emph{superlinearly convergent} w.r.t. \(\seq{y^k}\) if
	\begin{equation}\label{eq:SuperlinearDir}
		\lim_{k\to\infty}{
			\frac{\smash{\|y^k+d^k-y_\star\|}}{\|y^k-y_\star\|} = 0
		}.
	\end{equation}
\end{defin}

When \(y_\star\) is a strong minimizer, by \cite[Cor. 3.6]{drusvyatskiy2016error} the error bound \eqref{eq:ErrorBound} holds for some \(\beta,\nu > 0\) and \(Y_\star = \set{y_\star}\).
This, by \cref{it:subseqConv}, implies \(y^k\to y_\star\).
Therefore we have the following result.

\begin{thm}\label{thm:SuperlinearConv}
	Suppose that \(f\) and \(g\) satisfy \Cref{ass:fg2}, and that \eqref{eq:DualProblem} has a (unique) strong minimizer \(y_\star\).
	If \eqref{eq:SuperlinearDir} holds in \Cref{alg:NAMA}, then
	\begin{enumerate}
		\item\label{it:SuperlinearStep} the stepsize \(\tau_k = 1\) for all \(k\) sufficiently large,
		\item\label{it:SuperlinearCost} the cost \(\psi(y^k)\to \inf\psi\) Q-superlinearly,
		\item\label{it:SuperlinearDual} the dual iterates \(y^k\to y_\star\) Q-superlinearly,
		\item\label{it:SuperlinearPrimal} the primal iterates \(x^k\to x_\star\) R-superlinearly.
	\end{enumerate}
\begin{proof}
We know from \cref{it:psigammaTwiceDiff,thm:HessDef+} that \(\psi_\gamma\) is twice differentiable with symmetric and positive definite Hessian \(H_\star=\nabla^2\psi_\gamma(y_\star)\).
We can expand \(\psi_\gamma\) around \(y_\star\) and obtain
\ifieee
	\begin{align*}
	&
	\frac{\psi_\gamma(y^k+d^k) - \inf\psi}{\psi_\gamma(y^k) - \inf\psi}
	\\
	{={}} &
		\frac{\innprod{H_\star (y^k+d^k-y_\star)}{y^k+d^k-y_\star}+o(\|y^k+d^k-y_\star\|^2)}{\innprod{H_\star(y^k-y_\star)}{y^k-y_\star}+o(\|y^k-y_\star\|^2)}
	\\
	{\leq{}} &
		\frac{
			\|H_\star\|
			\left(
				\frac{
					\|y^k+d^k-y_\star\|
				}{
					\|y^k-y_\star\|
				}
			\right)^2
			{}+{}
			\left(
				\frac{
					o(\|y^k+d^k-y_\star\|)
				}{
					\|y^k-y_\star\|
				}
			\right)^2
		}{
			\lambda_{\min}(H_\star)
			{}+{}
			\left(
				\frac{
					o(\|y^k-y_\star\|)
				}{
					\|y^k-y_\star\|
				}
			\right)^2
		}
	\end{align*}
\else
	\begin{align*}
		\frac{\psi_\gamma(y^k+d^k) - \inf\psi}{\psi_\gamma(y^k) - \inf\psi}
	{}={} &
		\frac{\innprod{H_\star (y^k+d^k-y_\star)}{y^k+d^k-y_\star}+o(\|y^k+d^k-y_\star\|^2)}{\innprod{H_\star(y^k-y_\star)}{y^k-y_\star}+o(\|y^k-y_\star\|^2)}
	\\
	{}\leq{} &
		\frac{
			\|H_\star\|
			\left(
				\frac{
					\|y^k+d^k-y_\star\|
				}{
					\|y^k-y_\star\|
				}
			\right)^2
			{}+{}
			\left(
				\frac{
					o(\|y^k+d^k-y_\star\|)
				}{
					\|y^k-y_\star\|
				}
			\right)^2
		}{
			\lambda_{\min}(H_\star)
			{}+{}
			\left(
				\frac{
					o(\|y^k-y_\star\|)
				}{
					\|y^k-y_\star\|
				}
			\right)^2
		}
	\end{align*}
\fi
which vanishes for \(k\to\infty\).
In particular, eventually \(\psi_\gamma(y^k+d^k)\leq\psi_\gamma(y^k)\) will always hold, proving \ref{it:SuperlinearStep}.
In turn, since eventually \(\tilde y^k=y^k+\tau_kd^k=y^k+d^k\), using \cref{it:AMEgeq} and \eqref{eq:chain2} we have
\[
	\frac{\psi(y^{k+1}) - \inf\psi}{\psi(y^k) - \inf\psi}
{}\leq{}
	\frac{\psi_\gamma(\tilde y^k) - \inf\psi}{\psi_\gamma(y^k) - \inf\psi}
{}\to{}
	0,
\]
which proves \ref{it:SuperlinearCost}.
Moreover, \eqref{eq:SuperlinearDir} reads
\begin{equation}\label{eq:SuperlinearTilde}
	\|\tilde y^k-y_\star\|/\|y^k-y_\star\| \to 0.
\end{equation}
Now, using nonexpansiveness of \(T_\gamma\) (cf. the proof of \cite[Thm. 25.8]{bauschke2011convex}) one has
\begin{align*}
	\|y^{k+1} {-{}} y_\star\|
{}={}
	\|T_\gamma(\tilde y^k) - T_\gamma(y_\star)\|
{}\leq{}
	\|\tilde y^k - y_\star\|
\end{align*}
which, with \eqref{eq:SuperlinearTilde}, proves \ref{it:SuperlinearDual}.
\ref{it:SuperlinearPrimal} follows from \ref{it:SuperlinearCost} and \cref{lem:dual2primal}.
\end{proof}
\end{thm}

		When quasi-Newton directions are computed as in \eqref{eq:QNdir}, superlinear convergence holds provided that the sequence of matrices \(\seq{B_k}\) satisfies the Dennis-Mor\'e condition given in the following result. Such condition is satisfied for example by the modified Broyden method \eqref{eq:ModifiedBroyden} under standard assumptions of \emph{calm semidifferentiability} of \(R_\gamma\), see \cite[Thm. 6.8]{themelis2016supermann}.

\begin{thm}[Dennis-Mor\'e condition]\label{thm:DennisMore}
	Suppose that \(f\) and \(g\) strictly satisfy \Cref{ass:fg2}, and that \eqref{eq:DualProblem} has a (unique) strong minimizer \(y_\star\).
	If \(\seq{d^k}\) is selected according to \eqref{eq:QNdir}, with
	\begin{equation}\label{eq:DennisMore}
		\lim_{k\to\infty}{
			\frac{\smash{\|(B_k - \jac{R_\gamma}{y_\star}) d^k\|}}{\|d^k\|} = 0
		},
	\end{equation}
	then \(\seq{d^k}\) is superlinearly convergent with respect to \(\seq{y^k}\).
	In particular, the conclusions of \Cref{thm:SuperlinearConv} hold.
	\begin{proof}
	From \cref{it:Rdiff,thm:JRnonsing} we know that \(R_\gamma\) is strictly differentiable, with nonsingular Jacobian \(J_\star = \jac{R_\gamma}{y_\star}\).
	Let us denote \(r^k = z^k - Ax^k = R_\gamma(y^k)\) for simplicity.
	By using \eqref{eq:QNdir} and \eqref{eq:DennisMore}, and by applying the reverse triangle inequality we obtain
	\begin{equation*}
			0 {}\leftarrow{} \frac{\|r^k - J_\star d^k\|}{\|d^k\|}
		{}\geq{}
			\frac{\|J_\star B_k^{-1}r^k\|}{\|d^k\|} - \frac{\|r^k\|}{\|d^k\|}
		{}\geq{}
			\alpha - \frac{\|r^k\|}{\|d^k\|},
	\end{equation*}
	where \(\alpha = \sqrt{\lambda_{\min}(\trans{J_\star}J_\star)} > 0\) since  \(J_\star\) is nonsingular. Therefore,
	\[
		\liminf_{k\to \infty}{
			\nicefrac{\|r^k\|}{\|d^k\|}
		}
	{}\geq{}
		\alpha
	\]
	and as a consequence \(\|d^k\| \leq (2/\alpha)\|r^k\|\) for all \(k\) sufficiently large.
	Since \(r^k \to 0\) by \cref{it:subseqConv}, then \(d^k \to 0\).
	We have
		\[
			\mathtight
			0 {}\leftarrow{} \frac{r^k - J_\star d^k}{\|d^k\|}
		{}={}
			\frac{r^k + J_\star d^k {-{}} R_\gamma(y^k+d^k)}{\|d^k\|} + \frac{R_\gamma(y^k+d^k)}{\|d^k\|}.
		\]
	The first summand in the above equation tends to zero because of strict differentiability of \(R_\gamma\) at \(y_\star\), therefore
	\[ R_\gamma(y^k+d^k)/\|d^k\| \to 0. \]
	By nonsingularity of \(J_\star\) then \(\|R_\gamma(y)\| \geq \alpha\|y-y_\star\|\) for all \(y\) sufficiently close to \(y_\star\), and
	since \(y^k+d^k\to y_\star\) we have
	\ifieee
		\begin{align*}
			0
		{}\leftarrow{}
			\frac{R_\gamma(y^k+d^k)}{\|d^k\|}
		{}\geq{} &
			\frac{\alpha\|y^k + d^k - y_\star\|}{\|d^k\|}
		\\
		{}\geq{} &
			\frac{\alpha\|y^k + d^k - y_\star\|}{\|y+d^k-y_\star\| + \|y^k-y_\star\|}.
		\end{align*}
	\else
		\[
			0
		{}\leftarrow{}
			\frac{R_\gamma(y^k+d^k)}{\|d^k\|}
		{}\geq{}
			\frac{\alpha\|y^k + d^k - y_\star\|}{\|d^k\|}
		{}\geq{}
			\frac{\alpha\|y^k + d^k - y_\star\|}{\|y+d^k-y_\star\| + \|y^k-y_\star\|}.
		\]
	\fi
	This implies \(\|y^k + d^k - y_\star\|/\|y^k-y_\star\| \to 0\),
	\ie \(\seq{d^k}\) is superlinearly convergent with respect to \(\seq{y^k}\).
	\end{proof}
\end{thm}


	\section{Simulations}
		\label{sec:Simulations}
		We now present numerical results obtained with the proposed algorithm.
The scripts reproducing the results in this section are available online.\footnote{\url{https://github.com/kul-forbes/NAMA-experiments}}
In \algacronym\ we used \(\beta = 0.5\) and \(\tau_{\min} = 10^{-3}\) (see \Cref{rem:IterationComplexity}). Furthermore,
in all experiments we computed directions \(\seq{d^k}\) according to the L-BFGS method, with memory \(20\), which is able to scale with the problem dimension much better then full quasi-Newton update formulas.
All experiments were performed using MATLAB 2016b (v9.1.0) on a MacBook Pro running macOS 10.12, with an Intel Core i5 CPU (2.7 GHz) and 8 GB of memory.

		\subsection{Linear MPC}
			We consider finite horizon, discrete time, linear optimal control problems of the form
\begin{subequations}\label{eq:MPCprob}
\begin{align}
	\minimize_{\substack{x_0,\ldots,x_N \\ u_0,\ldots,u_{N-1}}}\ &{} \sum_{i=0}^{N-1}\ell_i(x_i,u_i)+\ell_N(x_N)\\
	\stt\ &{} x_0=x_{\textrm{init}},\\
		  &{} x_{i+1}=\Phi_i x_i + \Gamma_i u_i + c_i,\ i=0,\ldots,N-1, \label{eq:MPCdyn}
\shortintertext{where \(x_0,\ldots,x_N\in\R^{n_x}\) and \(u_0,\ldots,u_{N-1}\in\R^{n_u}\), and}
	& \ell_{i}(x,u) = q_i(x,u)+g_i(L_i(x,u)),\label{eq:MPCcost}\\
	& \ell_N(x) = q_N(x)+g_N(L_N x).
\end{align}\end{subequations}
Here the $q_i$ are strongly convex (typically quadratic), the $g_i$ are proper, closed, convex functions, while the \(L_i\) are linear mappings, for \(i=0,\ldots,N\). For example, with a convex set \(C\), one can set
\begin{align}
	g_i(\cdot) &= \indicator_{C}(\cdot) \hspace{28mm} \mbox{(hard constraints)} \nonumber\\
	g_i(\cdot) &= \alpha\dist_{C}(\cdot),\quad\alpha>0, \hspace{5mm} \mbox{(soft constraints)} \nonumber\\
\shortintertext{Set $C$ here is typically the nonpositive orthant or a box, but can be any other convex set onto which one can efficiently project. When \(C = [a_1, b_1]\times\ldots\times[a_d, b_d]\) is a \(d\)-dimensional box, then one can alternatively model soft constraints as}
	g_i(z) &{}= \textstyle\sum_{j=1}^d \alpha_j \bigl|z_j - \max\set{a_j, \min\set{b_j, z_j}}\bigr|. \label{eq:MPCsoftconstr}
\end{align}

Problem \eqref{eq:MPCprob} takes the form \eqref{eq:Problem} by reformulating it as follows (see also \cite{patrinos2014accelerated,giselsson2015metric,stathopoulos2016splitting}). Denote the full sequence of states and inputs as \(\bar x = (x_0,u_0,x_1,u_1,\ldots, x_N)\), and let
\[S(p) = \set{\bar x}[x_i = \Phi_i x_i + \Gamma_i u_i, x_0 = p]\]
be the affine subspace of feasible trajectories of the system having initial state \(p\). Then in \eqref{eq:Problem}
\begin{align*}
	f(\bar x) &{}={} \textstyle\sum_{i=0}^{N-1} q_i(x_i,u_i) + q_N(x_N) + \delta_{S(x_{\textrm{init}})}(\bar x), \\
	g(\bar z) &{}={} \textstyle\sum_{i=0}^{N} g_i(z_i),\quad A {}={} \diag(L_0,\ldots,L_N).
\end{align*}

Let us further denote by \(\bar y = (y_0,\ldots,y_N)\) the dual variable associated with the above problem.
In this case, in the alternating minimization step \ref{step:AMLS1_altmin} of \algacronym{}, the iterate \(\bar x^k\) is obtained by solving
\begin{align*}
	\minimize\ & \textstyle\sum_{i=0}^{N-1} q_i(x_i,u_i) + \innprod{y_i^k}{L_i(x_i, u_i)} \\
	& + q_N(x_N) + \innprod{y_N^k}{L_N x_N}.\\
	\stt\ & x_{i+1}=\Phi_i x_i + \Gamma_i u_i + c_i,\ i=0,\ldots,N-1.
\end{align*}
This is an unconstrained LQR problem whose solution can be efficiently computed with a Riccati-like recursion procedure, in the typical case where \(q_0,\ldots,q_N\) are quadratic, see \cite[Alg.s 3, 4]{patrinos2014accelerated}. The expensive ``factor'' step only needs to be performed once, before the main loop of the algorithm takes place. At every iteration one needs to perform merely a forward-backward sweep and no matrix inversions are required. Furthermore
\begin{align*}
	\bar z_i^k {}={}& \prox_{\gamma^{-1}g_i}(\gamma^{-1}y_i^k + L_i(x_i^k,u_i^k)), \quad i=0,\ldots,N-1, \\
	\bar z_N^k {}={}& \prox_{\gamma^{-1}g_N}(\gamma^{-1}y_N^k + L_N(x_N^k)),
\end{align*}
which in the case of hard/soft constraints essentially consist of projections onto the constrained sets.

\subsubsection{Aircraft control}

We applied the proposed method to the AFTI-16 aircraft control problem \cite{bemporad1997nonlinear,giselsson2015metric} with \(n_x = 4\) states and \(n_u = 2\) inputs, for a sampling time \(T_s = 0.05\) seconds.
The objective is to drive the \emph{pitch angle} from \(0^\circ\) to \(10^\circ\), and then back to \(0^\circ\). We simulated the system for \(4\) seconds, at the sampling time \(T_s = 0.05\), using \(N=50\) and quadratic costs
\begin{align*}
	q_i(x, u) &= \tfrac12\|x-x_{\textrm{ref}}\|_Q^2 + \tfrac12\|u\|_R^2,\quad i=0,\ldots,N-1, \\
	q_N(x) &= \tfrac12\|x-x_{\textrm{ref}}\|_{Q_N}^2,
\end{align*}
where \(Q = \diag(10^{-4}, 10^2, 10^{-3}, 10^2)\), \(Q_N = 100\cdot Q\) and \(R = \diag(10^{-2}, 10^{-2})\). The reference was set \(x_{\textrm{ref}} = (0, 0, 0, 10)\) for the first \(2\) seconds, and \(x_{\textrm{ref}} = (0, 0, 0, 0)\) for the remaining \(2\) seconds. Furthermore, we imposed hard box constraints on the inputs, and soft box constraints \eqref{eq:MPCsoftconstr} on the states, with weights \(10^6\).
Since soft constraints can be formulated into a QP, by adding linearly penalized nonnegative slack variables, we also compared against standard QP solvers.

The dual problem has a condition number of \(10^8\). To improve the convergence of the algorithms we therefore considered scaling the dual variables according to the \emph{Jacobi scaling}, which consists of a diagonal change of variable (in the dual space) enforcing the (dual) Hessian to have diagonal elements equal to one (see also \cite{richter2013certification,giselsson2015metric} on the problem of preconditioning fast dual proximal gradient methods). Note that a diagonal change of variable in the dual space simply corresponds to a scaling of the equality constraints, when the problem is equivalently formulated as \eqref{eq:AltProblem}.

We compared \algacronym\ against fast AMA \cite{pu2017complexity}, which is also known as GPAD \cite{patrinos2014accelerated} in this context, qpOASES v3.2.0 \cite{ferreau2014qpoases} and the commercial QP solver MOSEK v7.1.
We also compared against the cone solvers ECOS v2.0.4 \cite{domahidi2013ecos}, SDPT3 v4.0 \cite{toh1999sdpt3} and SeDuMi v1.34 \cite{sturm1999using}, all accessed through CVX v2.1 in MATLAB: note that the CPU time for these methods does not include the problem parsing and preprocessing by CVX, but only considers the actual running time of the solvers.
The results of the simulations are reported in \Cref{tbl:AFTI16}. As termination criterion for \algacronym\ and GPAD we used \(\|R_\gamma(y^k)\|_\infty \leq \epsilon_{\textrm{tol}} = 10^{-4}\). We also report the (average and maximum) number of \(x\)- and \(z\)-minimization steps performed by NAMA: due to the structure of \(f\), the \(x\)-update is a linear mapping, and consequently we can save its computation during the backtracking line-search. GPAD, in contrast, performs one alternating minimization per iteration.

Apparently, \algacronym\ greatly improves the convergence performance with respect to GPAD.
When the problem is prescaled, our method performs favorably also with respect to the other QP and cone solvers considered.
One must keep in mind that \algacronym\ was executed using a generic, high-level MATLAB implementation. As computation times become smaller and smaller, overheads due to the runtime environment get more and more relevant in the total CPU time.
A tailored, low-level implementation of the same algorithm could significantly decrease the CPU times shown in \Cref{tbl:AFTI16}: this is also reported in \cite{giselsson2015metric}, where a speedup of more than a factor \(20\) is observed using C code generation.

\begin{table*}\centering
{\small
\rowcolors{2}{gray!0}{gray!10}
\begin{tabular}{|ll|rr|rr|rr|rr|}
\rowcolor{gray!25}
\hline
			&						& \multicolumn{2}{c|}{Iterations}	& \multicolumn{2}{c|}{\(x\)-updates}	& \multicolumn{2}{c|}{\(z\)-updates}	& \multicolumn{2}{c|}{CPU time (ms)} \\
\rowcolor{gray!25}
			&						& avg. & max.			& avg. & max.	 	& avg. & max.   & avg. & max. \\
\hline
GPAD 	& (no scaling) 			& 6408.2 & 118.3 k & - & - & - & - & 1645.7 & 23331.9 \\
\algacronym\ (L-BFGS, mem = 20) & (no scaling) 			& 66.0 & 748 & 134.2 & 1527 & 139.7 & 1565 & 36.5 & 464.6 \\
\hline
GPAD 	& (Jacobi scaling) 		& 104.8 & 491 & - & - & - & - & 21.0 & 96.7 \\
\textbf{\algacronym\ (L-BFGS, mem = 20)} & \textbf{(Jacobi scaling)}		& \textbf{9.7} & \textbf{42} & \textbf{18.7} & \textbf{85} & \textbf{18.8} & \textbf{88} & \textbf{4.9} & \textbf{21.3} \\
\hline
qpOASES 	&                		& 	 	& 	&  &  &  & 	 & 2362.7 & 2573.3 \\
qpOASES		& (warm-started)		& 		&		&  &  &  &  & 14.6 & 286.9 \\
MOSEK		& 					& 		&		&  &  &  &  & 207.4 & 539.4 \\
ECOS 	&                		& 	 	& 	&  &  &  & 	 & 23.6 & 37.6 \\
SDPT3 	&                		& 	 	& 	&  &  &  & 	 & 607.7 & 890.6 \\
SeDuMi 	&                		& 	 	& 		&  &  &  &  & 137.2 & 266.2 \\
\hline
\end{tabular}}

\caption{Aircraft control, performance of the algorithms in the case of the AFTI-16 problem, for \(T_s = 50\) ms and \(N=50\). GPAD and \algacronym\ were stopped as soon as \(\|R_\gamma(y^k)\|_\infty \leq \epsilon_{\textrm{tol}}=10^{-4}\). Since the problem is ill-conditioned, we also applied the methods by prescaling the dual problem.
The number of \(x\)- and \(z\)- updates of GPAD equals the number of iterations.
\algacronym\ was executed using a generic implementation in MATLAB, while the others QP and cone solvers considered are all implemented in C/C++.}\label{tbl:AFTI16}
\end{table*}

\subsubsection{Oscillating masses}

\begin{figure}
	\centering
	{{%
		\pgfkeys{/pgf/images/include external/.code={\includegraphics[width=0.4\textwidth]{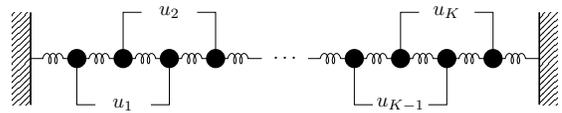}}}%
		\tikzsetnextfilename{chain}%
		\begin{tikzpicture}[every node/.style={outer sep=0pt},scale=0.8]
\tikzstyle{spring}=[decorate,decoration={coil, pre length=6, post length=6, segment length=3}]
\fill[pattern = north east lines] (0,-1) rectangle (-0.4,1);
\draw[thick] (0,-1) -- (0,1);

\def \springlength {1.0}
\def \handlespace {0.6}

\foreach \k in {1}
{

    \def \basex {((\k - 1)*4)*\springlength}
    \draw[spring] ({\basex},0) -- ({\basex+\springlength},0);
    \node[circle,fill=black] (a) at ({\basex+\springlength},0) {};
    \draw[spring] ({\basex+\springlength},0) -- ({\basex+2*\springlength},0);
    \node[circle,fill=black] (a) at ({\basex+2*\springlength},0) {};
    \draw[spring] ({\basex+2*\springlength},0) -- ({\basex+3*\springlength},0);
    \node[circle,fill=black] (a) at ({\basex+3*\springlength},0) {};
    \draw[spring] ({\basex+3*\springlength},0) -- ({\basex+4*\springlength},0);
    \node[circle,fill=black] (a) at ({\basex+4*\springlength},0) {};


    \pgfmathsetmacro\ucount{(\k-1)*2+1}

    \draw ({\basex+\springlength},0) -- ({\basex+\springlength},-1) -- ({\basex+2*\springlength-\handlespace},-1);
    \draw ({\basex+3*\springlength},0) -- ({\basex+3*\springlength},-1) -- ({\basex+2*\springlength+\handlespace},-1);
    \node at ({\basex+2*\springlength},-1) {\(\:u_1\:\)}; 

    \pgfmathsetmacro\ucount{(\k-1)*2+2}

    \draw ({\basex+2*\springlength},0) -- ({\basex+2*\springlength},+1) -- ({\basex+3*\springlength-\handlespace},+1);
    \draw ({\basex+4*\springlength},0) -- ({\basex+4*\springlength},+1) -- ({\basex+3*\springlength+\handlespace},+1);
    \node at ({\basex+3*\springlength},+1) {\(\:u_2\:\)}; 
}

\draw[spring] ({4*\springlength},0) -- ({5*\springlength},0);
\node at ({5.5*\springlength},0) {\(\cdots\)};

\foreach \k in {1}
{

    \def \basex {(\k*4+2)*\springlength}
    \draw[spring] ({\basex},0) -- ({\basex+\springlength},0);
    \node[circle,fill=black] (a) at ({\basex+\springlength},0) {};
    \draw[spring] ({\basex+\springlength},0) -- ({\basex+2*\springlength},0);
    \node[circle,fill=black] (a) at ({\basex+2*\springlength},0) {};
    \draw[spring] ({\basex+2*\springlength},0) -- ({\basex+3*\springlength},0);
    \node[circle,fill=black] (a) at ({\basex+3*\springlength},0) {};
    \draw[spring] ({\basex+3*\springlength},0) -- ({\basex+4*\springlength},0);
    \node[circle,fill=black] (a) at ({\basex+4*\springlength},0) {};


    \pgfmathsetmacro\ucount{(\k-1)*2+1}

    \draw ({\basex+\springlength},0) -- ({\basex+\springlength},-1) -- ({\basex+2*\springlength-\handlespace},-1);
    \draw ({\basex+3*\springlength},0) -- ({\basex+3*\springlength},-1) -- ({\basex+2*\springlength+\handlespace},-1);
    \node at ({\basex+2*\springlength},-1) {\(\:u_{K-1}\:\)}; 

    \pgfmathsetmacro\ucount{(\k-1)*2+2}

    \draw ({\basex+2*\springlength},0) -- ({\basex+2*\springlength},+1) -- ({\basex+3*\springlength-\handlespace},+1);
    \draw ({\basex+4*\springlength},0) -- ({\basex+4*\springlength},+1) -- ({\basex+3*\springlength+\handlespace},+1);
    \node at ({\basex+3*\springlength},+1) {\(\:u_{K}\:\)}; 
}

\def \basex {(2*4+2)*\springlength};

\draw[spring] ({\basex}, 0) -- ({\basex+\springlength}, 0);
\fill[pattern = north east lines] ({\basex+\springlength},-1) rectangle ({\basex+\springlength+0.4},1);
\draw[thick] ({\basex+\springlength},-1) -- ({\basex+\springlength},1);

%
%
%
%
%
%
%
%
%
\end{tikzpicture}%
	}}%
	\caption[Oscillating masses, schematic representation.]{Oscillating masses, schematic representation of the simulated system.
	}\label{fig:Masses}
\end{figure}

Next, we consider a chain of oscillating masses connected by springs, with both ends attached to walls. The chain is composed of \(2K\) bodies of unit mass, the springs have constant \(1\) and no damping, and the system is controlled through \(K\) actuators, each being a force acting on a pair of masses, as depicted in \Cref{fig:Masses}. Therefore \(n_x=4K\) (the states are the displacement from the rest position and velocity of each mass) and \(n_u=K\).
The inputs are constrained in \([-0.5,+0.5]\), while the position and velocity of each mass is constrained in \([-4,+4]\).

The continuous-time system was discretized with a sampling time \(T_s = 0.5\).
Like in the previous example, we considered quadratic costs with \(Q = Q_N = I_{n_x}\), \(R = I_{n_u}\) and hard constraints on state and input. Furthermore, we imposed a quadratic terminal constraint
\begin{equation}\label{eq:MPCterminal}\tfrac{1}{2}\innprod{Px_N}{x_N} \leq \delta,\end{equation}
where \(P\) solves the Riccati equation related to the discrete-time LQR problem.
Constraint \eqref{eq:MPCterminal} can be enforced by taking \(L_N\) in \eqref{eq:MPCprob} as the Cholesky factor of \(P\), so that \(\trans{L_N}L_N = P\), and \(g_N\) 
as the indicator of the Euclidean ball of radius \(\sqrt{\delta}\).
Parameter \(\delta\) is selected so as to ensure that no constraints are violated in such ellipsoidal set.

\begin{figure*}
	\centering
	{{%
		\pgfkeys{/pgf/images/include external/.code={\includegraphics[]{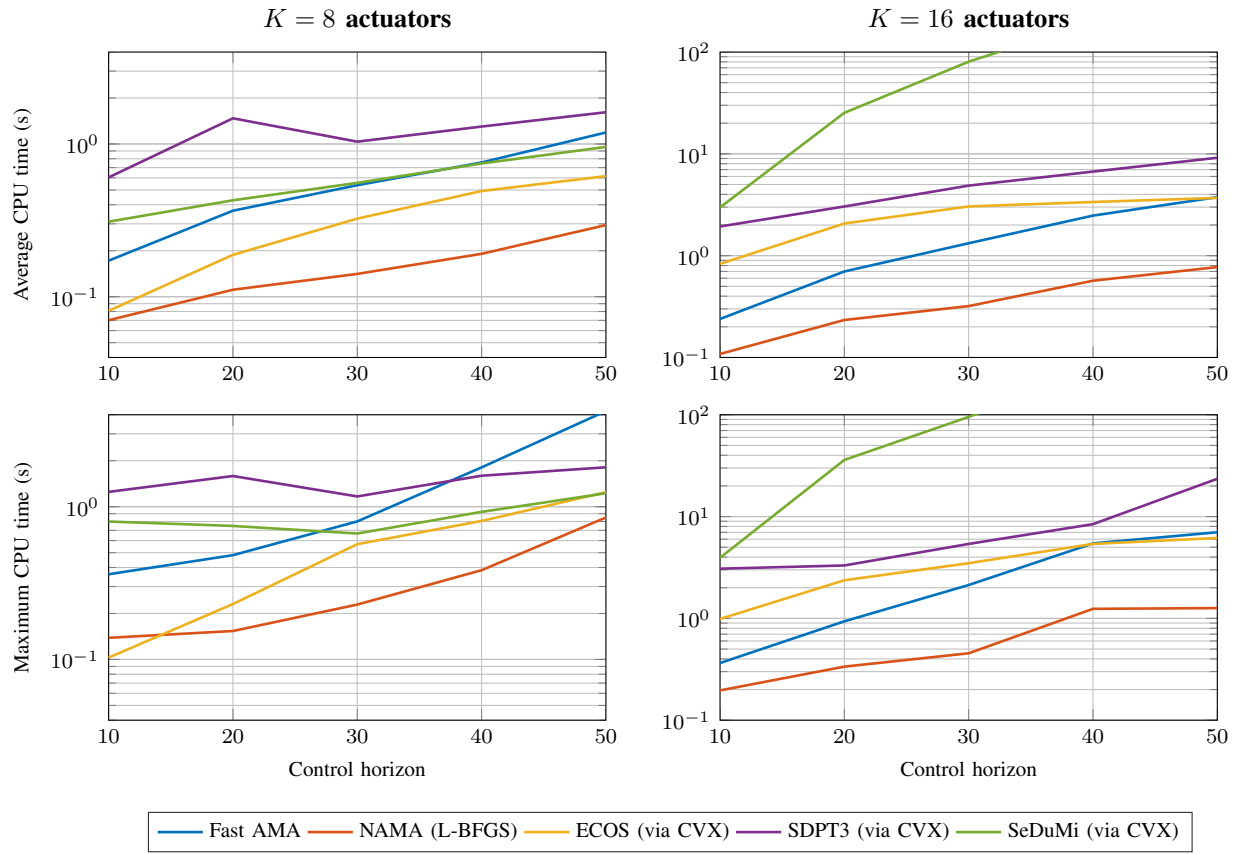}}}%
		\tikzsetnextfilename{masses_all}%
%
%
\definecolor{mycolor1}{rgb}{0.00000,0.44700,0.74100}%
\definecolor{mycolor2}{rgb}{0.85000,0.32500,0.09800}%
\definecolor{mycolor3}{rgb}{0.92900,0.69400,0.12500}%
\definecolor{mycolor4}{rgb}{0.49400,0.18400,0.55600}%
\definecolor{mycolor5}{rgb}{0.46600,0.67400,0.18800}%
\begin{tikzpicture}


\begin{axis}[%
width=2.603in,
height=1.600in,
at={(0.0in,0.0in)},
scale only axis,
xmin=10,
xmax=50,
xmajorgrids,
ymode=log,
ymin=0.04,
ymax=4,
yminorticks=true,
ylabel={Average CPU time (s)},
ymajorgrids,
yminorgrids,
axis background/.style={fill=white},
title style={font=\bfseries},
title={\(K=8\) actuators},
]
\addplot [color=mycolor1,solid,line width=1.0pt]
  table[row sep=crcr]{%
10	0.17218789092\\
20	0.36568405216\\
30	0.53706685456\\
40	0.75670134708\\
50	1.18909343792\\
};

\addplot [color=mycolor2,solid,line width=1.0pt]
  table[row sep=crcr]{%
10	0.07015270834\\
20	0.11109088422\\
30	0.14099680242\\
40	0.1908991365\\
50	0.2944189923\\
};

\addplot [color=mycolor3,solid,line width=1.0pt]
  table[row sep=crcr]{%
10	0.0808494815\\
20	0.18803111178\\
30	0.32444725404\\
40	0.4920617083\\
50	0.61581052024\\
};

\addplot [color=mycolor4,solid,line width=1.0pt]
  table[row sep=crcr]{%
10	0.60407715912\\
20	1.4724735567\\
30	1.03555906388\\
40	1.30023440908\\
50	1.61337119532\\
};

\addplot [color=mycolor5,solid,line width=1.0pt]
  table[row sep=crcr]{%
10	0.30986802676\\
20	0.4278456287\\
30	0.556173464\\
40	0.74546291272\\
50	0.95633066358\\
};

\end{axis}


\begin{axis}[%
width=2.603in,
height=1.600in,
at={(0.0in,-1.9in)},
scale only axis,
xmin=10,
xmax=50,
xlabel={Control horizon},
xmajorgrids,
ymode=log,
ymin=0.04,
ymax=4,
yminorticks=true,
ylabel={Maximum CPU time (s)},
ymajorgrids,
yminorgrids,
axis background/.style={fill=white},
]
\addplot [color=mycolor1,solid,line width=1.0pt]
  table[row sep=crcr]{%
10	0.361109715\\
20	0.481621875\\
30	0.801011958\\
40	1.808817828\\
50	4.227485807\\
};

\addplot [color=mycolor2,solid,line width=1.0pt]
  table[row sep=crcr]{%
10	0.138532844\\
20	0.153497935\\
30	0.228688196\\
40	0.383901042\\
50	0.849925307\\
};

\addplot [color=mycolor3,solid,line width=1.0pt]
  table[row sep=crcr]{%
10	0.102782544\\
20	0.230491504\\
30	0.568380803\\
40	0.806687013\\
50	1.244293179\\
};

\addplot [color=mycolor4,solid,line width=1.0pt]
  table[row sep=crcr]{%
10	1.250076603\\
20	1.588014395\\
30	1.167265807\\
40	1.595182224\\
50	1.811337881\\
};

\addplot [color=mycolor5,solid,line width=1.0pt]
  table[row sep=crcr]{%
10	0.799702143\\
20	0.748077432\\
30	0.66730216\\
40	0.924990083\\
50	1.224832753\\
};

\end{axis}


\begin{axis}[%
width=2.603in,
height=1.600in,
at={(3.2in,0.0in)},
scale only axis,
xmin=10,
xmax=50,
xmajorgrids,
ymode=log,
ymin=0.1,
ymax=100,
yminorticks=true,
ymajorgrids,
yminorgrids,
axis background/.style={fill=white},
title style={font=\bfseries},
title={\(K=16\) actuators},
]
\addplot [color=mycolor1,solid,line width=1.0pt]
  table[row sep=crcr]{%
10	0.23881885692\\
20	0.69934740704\\
30	1.32436439538\\
40	2.47437136486\\
50	3.75470422598\\
};

\addplot [color=mycolor2,solid,line width=1.0pt]
  table[row sep=crcr]{%
10	0.10828978842\\
20	0.23344533878\\
30	0.31976621136\\
40	0.56746490212\\
50	0.77273548824\\
};

\addplot [color=mycolor3,solid,line width=1.0pt]
  table[row sep=crcr]{%
10	0.83036378994\\
20	2.06909152476\\
30	3.04412290802\\
40	3.36550807742\\
50	3.692248494\\
};

\addplot [color=mycolor4,solid,line width=1.0pt]
  table[row sep=crcr]{%
10	1.93458965726\\
20	3.04112135144\\
30	4.8744867187\\
40	6.69654472298\\
50	9.1529801118\\
};

\addplot [color=mycolor5,solid,line width=1.0pt]
  table[row sep=crcr]{%
10	2.97531748292\\
20	25.2918726238\\
30	80.36843505492\\
40	200.06453938122\\
50	413.59278211554\\
};

\end{axis}


\begin{axis}[%
width=2.603in,
height=1.600in,
at={(3.2in,-1.9in)},
scale only axis,
xmin=10,
xmax=50,
xlabel={Control horizon},
xmajorgrids,
ymode=log,
ymin=0.1,
ymax=100,
yminorticks=true,
ymajorgrids,
yminorgrids,
axis background/.style={fill=white},
legend style={legend cell align=left,align=left,draw=white!15!black,at={(-0.1,-0.3)},anchor=north},
legend columns=5,
]
\addplot [color=mycolor1,solid,line width=1.0pt]
  table[row sep=crcr]{%
10	0.364362103\\
20	0.935794231\\
30	2.127214059\\
40	5.454526214\\
50	7.014129551\\
};
\addlegendentry{Fast AMA};

\addplot [color=mycolor2,solid,line width=1.0pt]
  table[row sep=crcr]{%
10	0.196182651\\
20	0.335944032\\
30	0.454267612\\
40	1.241576841\\
50	1.261734887\\
};
\addlegendentry{NAMA (L-BFGS)};

\addplot [color=mycolor3,solid,line width=1.0pt]
  table[row sep=crcr]{%
10	0.98744404\\
20	2.367935581\\
30	3.483973871\\
40	5.387731706\\
50	6.163720603\\
};
\addlegendentry{ECOS (via CVX)};

\addplot [color=mycolor4,solid,line width=1.0pt]
  table[row sep=crcr]{%
10	3.077851763\\
20	3.310074665\\
30	5.387058236\\
40	8.424509375\\
50	23.507160719\\
};
\addlegendentry{SDPT3 (via CVX)};

\addplot [color=mycolor5,solid,line width=1.0pt]
  table[row sep=crcr]{%
10	3.940593604\\
20	35.979146929\\
30	95.29951122\\
40	268.042432358\\
50	605.671043816\\
};
\addlegendentry{SeDuMi (via CVX)};

\end{axis}


\end{tikzpicture}%
%
	}}
	\caption[Oscillating masses, average and maximum CPU time.]{Oscillating masses, average and maximum CPU time (in seconds) for increasing prediction horizon and \(50\) randomly selected initial states. First column: \(K=8\) actuators. Second column: \(K=16\) actuators. Fast AMA and \algacronym\ were stopped as soon as \(\|R_\gamma(y^k)\|_\infty \leq \epsilon_{\textrm{tol}}=10^{-4}\).}\label{fig:ResultsMasses}
\end{figure*}

We simulated different scenarios, each with a different prediction horizon
\(N \in \set{10, 20, \ldots, 50}\), with \(K=8,16\).
For each scenario we selected \(50\) random initial states \(x_{\textrm{init}}\) by solving random feasibility problems (\eg with a cone solver) so as to ensure that a feasible trajectory starting from \(x_{\textrm{init}}\) exists.
Every algorithm was executed with the same set of initial conditions.
The results of this experiment are shown in \Cref{fig:ResultsMasses}. In addition to fast AMA, we compared \algacronym\ against ECOS, SDPT3 and SeDuMi, all accessed through CVX in MATLAB. \algacronym\ compares considerably well with all the other methods in this example, and in particular outperforms fast AMA, both on average and in the worst case.

	\section{Conclusions}
		In this work we presented \algacronym, a line-search method for minimizing the sum of two convex functions, one of which is assumed to be strongly convex, while the other is composed with a linear transformation.
The method is an extension of the classical alternating minimization algorithm (AMA), performing an additional line-search step over the \emph{alternating minimization envelope} associated with the problem.
By appropriately selecting the line-search directions, for example according to quasi-Newton methods for solving the optimality conditions \(R_\gamma(y) = 0\), we have shown that the algorithm converges superlinearly provided that ordinary second-order sufficiency conditions hold for the envelope function at the (unique) dual solution.
At the same time, the algorithm possesses the same global sublinear and local linear convergence rates as AMA.
Numerical experiments with the proposed method on linear MPC problems suggest that \algacronym\ is able to significantly speed up the convergence of AMA, comparing favorably against its accelerated variant and other state-of-the-art solvers even when limited-memory methods, such as L-BFGS, are used to compute the search directions.


	\begin{appendix}
		\ifieee\relax\else
			\section{Additional results}
		\fi
		\begin{lem}\label{lem:TwoPointsIneq}
	Let \(y, w\in\R^m\) and \(\gamma > 0\).
	Then,
	\begin{equation}\begin{aligned}
		\psi(w) &{}\geq{} \psi_\gamma(y)+\tfrac{\gamma}{2}\|Ax(y)-z_\gamma(y)\|^2 \\
		&\phantom{{}\geq{}} + \innprod{z_\gamma(y) - Ax(y)}{w-y}.
	\end{aligned}\end{equation}
\begin{proof}
	By \eqref{eq:FenchelIneq} we have
	\begin{align*}
			f(x(y))+\conj{f}(-\trans{A} w)
		&{}\geq{}
			-\innprod{ Ax(y)}{w},
		\\
			g(z_\gamma(y))+\conj{g}(w)
		&{}\geq{}
			\innprod{z_\gamma(y)}{w}.
	\end{align*}
	By summing the two inequalities and using the definition of \(\psi_\gamma\), after manipulations one obtains the result.
\end{proof}
\end{lem}

\begin{lem}\label{lem:dual2primal}
For all \(y\in\R^m\) it holds
\[
	\tfrac{\mu_f}{2}\|x(y)-x_\star\|^2 \leq \psi(y) - \inf\psi.
\]
\begin{proof}
From the optimality condition of the problem defining \(x(y)\), one obtains \(-\trans{A}y\in\partial f(x(y))\). Then, by strong convexity of \(f\) one gets
\[f(x(y)) - \innprod{\trans{A}y}{x_\star - x(y)} + \tfrac{\mu_f}{2}\|x(y) - x_\star\|^2 \leq f(x_\star). \]
By using \eqref{eq:ConjSubgr_f} in the above inequality we obtain
\[
	\tfrac{\mu_f}{2}\|x(y) - x_\star\|^2 - \innprod{Ax_\star}{y} {}\leq{} f(x_\star) + \conj{f}(-\trans{A}y),
\]
By using \eqref{eq:FenchelIneq} on \(g\) we have instead
\[
	\innprod{Ax_\star}{y} {}\leq{} g(Ax_\star) + \conj{g}(y).
\]
By summing the last two inequalities one obtains
\[
	\tfrac{\mu_f}{2}\|x(y) - x_\star\|^2 {}\leq{} f(x_\star) + g(Ax_\star) + \psi(y),
\]
and the claimed bound follows by strong duality.
\end{proof}
\end{lem}

\begin{lem}\label{lem:linearReg}
Suppose that the following hold for \eqref{eq:Problem}:
\begin{enumerate}
	\item\label{it:LemStrictFeas} \(A\relint(\dom f) \cap \relint(\dom g) \neq \emptyset \) (strict feasibility);
	\item\label{it:LemStrictComp} \(0\in\relint\partial(f+g\circ A)(x_\star)\) (strict complementarity).
\end{enumerate}
Then for any compact set \(U\) there is \(\kappa > 0\) such that
	\ifieee
		\[
			\dist(y, Y_\star)
		{}\leq{}
			\kappa
			\left[
				\dist(-\trans{A}y,\partial f(x_\star))
				{}+{}
				\dist(y,\partial g(Ax_\star))
			\right]
		\]
		holds for all \(y\in U\).
	\else
		\[\dist(y, Y_\star) \leq \kappa\left[ \dist(-\trans{A}y,\partial f(x_\star)) + \dist(y,\partial g(Ax_\star))\right],\quad \forall y\in U.\]
	\fi
\begin{proof}
	From \cref{it:LemStrictComp} it follows that
	\begin{align}
		0 {}\in{}& \relint\left[\partial f(x_\star) + \trans{A}\partial g(Ax_\star) \right] \nonumber \\
		 {}={}& \relint\partial f(x_\star) + \trans{A}\relint \partial g(Ax_\star).\label{eq:StrictInclusion}
	\end{align}
	In fact, the first inclusion is due to \cite[Thm 23.9]{rockafellar1997convex} in light of \cref{it:LemStrictFeas}, and the equality is due to \cite[Thm. 6.6]{rockafellar1997convex}.
	Consider \(W = \set{w}[-\trans{A}w \in \partial f(x_\star)] \subseteq \R^m\). From \eqref{eq:FirstOrderNecessary},
	\[
		Y_\star = W \cap \partial g(Ax_\star).
	\]
	Furthermore, using \eqref{eq:StrictInclusion} we obtain
	\[
		 \emptyset \neq \set{w}[-\trans{A}w \in \relint \partial f(x_\star)] = \relint W,
	\]
	where the equality is due to \cite[Thm. 6.7]{rockafellar1997convex},
	and the fact that \(\relint W \cap \relint \partial g(Ax_\star) \neq \emptyset \).
	By \cite[Cor. 5]{bauschke1999strong} then, we conclude that \(W\) and \(\partial g(Ax_\star)\) are boundedly linearly regular: for any compact set \(U\) there is \(\alpha>0\) such that for all \(y\in U\)
	\begin{equation}
		\dist(y, Y_\star) \leq \alpha\bigl[ \dist(y, W) + \dist(y,\partial g(Ax_\star))\bigr]. \label{eq:BLR1}
	\end{equation}

	Similarly, \eqref{eq:StrictInclusion} implies with \cite[Cor. 5]{bauschke1999strong} that the sets
	\(
		L = \set{(w, -\trans{A}w)}[w \in \R^m]
	\)
	and
	\(
		M = \R^m \times \partial f(x_\star)
	\)
	are boundedly linearly regular. Observe that
	\[ L\cap M = \set{(w, -\trans{A}w)}[-\trans{A}w \in \partial f(x_\star)]. \]
	Therefore, there is \(\beta>0\) such that for all \(y\in U\)
	\begin{align*}
		\dist(y, W) &{}\leq{} \dist((y,-\trans{A}y), L\cap M) \\
	&{}\leq{}
		\beta[\dist((y,-\trans{A}y), L)
		{{}+}
		\dist((y,-\trans{A}y), M)]
	\\
	&{}={} \beta\dist(-\trans{A}y, \partial f(x_\star)),
	\end{align*}
	where the second inequality is due to bounded linear regularity of \(L\) and \(M\), while the equality holds since \((y,-\trans{A}y)\in L\) and \(\dist((y,-\trans{A}y), M) = \dist(-\trans{A}y, \partial f(x_\star))\) for any \(y\).
	Using the above inequality in \eqref{eq:BLR1} yields the result.
\end{proof}
\end{lem}

\begin{lem}[Twice differentiability of \(\conj f\)]\label{lem:psi1TwiceDiff}
Suppose that \(f\) satisfies \Cref{ass:f2} for the primal-dual solution \((x_\star,y_\star)\).
Then \(\conj f\) is of class \(\cont^2\) around \(y_\star\), with
\[
	\nabla^2 \conj f(y_\star)
{}={}
	\geninv{H_f}.
\]
\begin{proof}
From \cite[Thm. 13.21]{rockafellar2011variational} we know that \(\conj f\) is twice epi-differentiable at \(v\) for \(x\in\partial\conj f(v)\) iff \(f\) is twice epi-differentiable at \(x\) for \(v\), with the relation
\begin{equation}\label{eq:twiceepiconj}
	\twiceepi[x]{\conj f}{v} = \conj{\left[\twiceepi[v]{f}{x}\right]}.
\end{equation}
The cited proof trivially extends to strict twice differentiability, and in fact \(\conj f\) turns out to be \emph{strictly} twice epi-differentiable at \(x_\star\).
Since \(\range(H_f)+S_f^\bot = \R^n\), by applying \eqref{eq:twiceepiconj} to \eqref{eq:fGenQuad} and conjugating \(\twiceepi[-\trans{A} y_\star]{f}{x_\star}\) by means of \cite[Prop. E.3.2.1]{hiriart2001fundamentals} we obtain that function \(\conj f\) has purely quadratic second epi-derivative (as opposed to generalized quadratic)
\[
	\twiceepi[x_\star]{\conj f}{-\trans{A} y_\star}[w]
{}={}
	\innprod{\smash{\geninv{(\proj_{S_f}H_f\proj_{S_f})}w}}{w}
{}\overrel{\eqref{eq:H}}{}
	\innprod{\smash{\geninv H_f} w}{w}
\]
which is everywhere finite in particular.
The proof now follows from \cite[Cor. 4.7]{poliquin1996generalized}.
\end{proof}
\end{lem}
With similar reasonings, the following result easily follows.
\begin{lem}[Twice epi-differentiability of \(\conj g\)]\label{lem:TwiceEpiDuality}
Suppose that \(g\) (strictly) satisfies \Cref{ass:g2} for a primal-dual solution \((x_\star,y_\star)\).
Then \(\conj g\) is (strictly) twice epi-dif\-fer\-en\-tiable at \(y_\star\) for \(A x_\star\).
More precisely, letting \(\bar S=S_g^\bot + \range(H_g)\),%
\begin{equation}\label{eq:EqSecondEpi}
	\twiceepi[A x_\star]{\conj g}{y_\star}
{}={}
	\conj{\left[\twiceepi[y_\star]{g}{A x_\star}\right]}
{}={}
	\innprod{\geninv{H_g}{}\cdot{}}{{}\cdot{}} + \delta_{\bar S}.
\end{equation}
\end{lem}

	\end{appendix}

	\section*{Acknowledgment}
	The authors would like thank Dmitriy Drusvyatskiy for his contribution to the proof of \Cref{lem:linearReg}.


	\bibliography{TeX/Bibliography.bib}

\begin{thebibliography}{10}
\providecommand{\url}[1]{#1}
\csname url@rmstyle\endcsname
\providecommand{\newblock}{\relax}
\providecommand{\bibinfo}[2]{#2}
\providecommand\BIBentrySTDinterwordspacing{\spaceskip=0pt\relax}
\providecommand\BIBentryALTinterwordstretchfactor{4}
\providecommand\BIBentryALTinterwordspacing{\spaceskip=\fontdimen2\font plus
\BIBentryALTinterwordstretchfactor\fontdimen3\font minus
  \fontdimen4\font\relax}
\providecommand\BIBforeignlanguage[2]{{%
\expandafter\ifx\csname l@#1\endcsname\relax
\typeout{** WARNING: IEEEtran.bst: No hyphenation pattern has been}%
\typeout{** loaded for the language `#1'. Using the pattern for}%
\typeout{** the default language instead.}%
\else
\language=\csname l@#1\endcsname
\fi
#2}}

\bibitem{stathopoulos2016splitting}
G.~Stathopoulos, H.~Shukla, A.~Szucs, Y.~Pu, and C.~N. Jones, ``Operator
  splitting methods in control,'' \emph{Foundations and Trends in Systems and
  Control}, vol.~3, no.~3, pp. 249--362, 2016.

\bibitem{fazel2013hankel}
M.~Fazel, T.~K. Pong, D.~Sun, and P.~Tseng, ``{H}ankel matrix rank minimization
  with applications to system identification and realization,'' \emph{SIAM
  Journal on Matrix Analysis and Applications}, vol.~34, no.~3, pp. 946--977,
  2013.

\bibitem{boyd2011distributed}
S.~Boyd, N.~Parikh, E.~Chu, B.~Peleato, and J.~Eckstein, ``Distributed
  optimization and statistical learning via the alternating direction method of
  multipliers,'' \emph{Foundations and Trends in Machine Learning}, vol.~3,
  no.~1, p. 1–122, 2011.

\bibitem{parikh2014proximal}
N.~Parikh and S.~Boyd, ``Proximal algorithms,'' \emph{Foundations and Trends in
  Optimization}, vol.~1, no.~3, pp. 127--239, 2014.

\bibitem{tseng1991applications}
P.~Tseng, ``Applications of a splitting algorithm to decomposition in convex
  programming and variational inequalities,'' \emph{{SIAM} Journal on Control
  and Optimization}, vol.~29, no.~1, pp. 119--138, 1991.

\bibitem{lions1979splitting}
P.-L. Lions and B.~Mercier, ``Splitting algorithms for the sum of two nonlinear
  operators,'' \emph{SIAM Journal on Numerical Analysis}, vol.~16, no.~6, pp.
  964--979, 1979.

\bibitem{beck2009fast}
A.~Beck and M.~Teboulle, ``A fast iterative shrinkage-thresholding algorithm
  for linear inverse problems,'' \emph{{SIAM} Journal on Imaging Sciences},
  vol.~2, no.~1, pp. 183--202, 2009.

\bibitem{nesterov2013gradient}
Y.~Nesterov, ``\BIBforeignlanguage{English}{Gradient methods for minimizing
  composite functions},'' \emph{\BIBforeignlanguage{English}{Mathematical
  Programming}}, vol. 140, no.~1, pp. 125--161, 2013.

\bibitem{beck2014fast}
A.~Beck and M.~Teboulle, ``\BIBforeignlanguage{English}{A fast dual proximal
  gradient algorithm for convex minimization and applications},''
  \emph{\BIBforeignlanguage{English}{Operations Research Letters}}, vol.~42,
  no.~1, pp. 1--6, 2014.

\bibitem{patrinos2013proximal}
P.~Patrinos and A.~Bemporad, ``Proximal {N}ewton methods for convex composite
  optimization,'' in \emph{IEEE Conference on Decision and Control}, 2013, pp.
  2358--2363.

\bibitem{stella2017forward}
L.~Stella, A.~Themelis, and P.~Patrinos, ``Forward-backward quasi-{N}ewton
  methods for nonsmooth optimization problems,'' \emph{Computational
  Optimization and Applications}, vol.~67, no.~3, pp. 443--487, 2017.

\bibitem{themelis2016forward}
A.~Themelis, L.~Stella, and P.~Patrinos, ``Forward-backward envelope for the
  sum of two nonconvex functions: Further properties and nonmonotone
  line-search algorithms,'' \emph{arXiv preprint arXiv:1606.06256}, 2016.

\bibitem{liu2016further}
T.~Liu and T.~K. Pong, ``Further properties of the forward--backward envelope
  with applications to difference-of-convex programming,'' \emph{Computational
  Optimization and Applications}, vol.~67, no.~3, pp. 489--520, 2017.

\bibitem{sampathirao2017proximal}
A.~K. Sampathirao, P.~Sopasakis, A.~Bemporad, and P.~Patrinos, ``Proximal
  limited-memory quasi-{N}ewton methods for scenario-based stochastic optimal
  control,'' \emph{To appear in Proceedings of the 20th IFAC Congress}, 2017.

\bibitem{patrinos2014douglas}
P.~Patrinos, L.~Stella, and A.~Bemporad, ``Douglas-{R}achford splitting:
  Complexity estimates and accelerated variants,'' in \emph{53rd IEEE
  Conference on Decision and Control}, 2014, pp. 4234--4239.

\bibitem{themelis2017douglas}
A.~Themelis, L.~Stella, and P.~Patrinos, ``Douglas--{R}achford splitting and
  {ADMM} for nonconvex optimization: new convergence results and accelerated
  versions,'' \emph{arXiv preprint arXiv:1709.05747}, 2017.

\bibitem{rockafellar1997convex}
R.~T. Rockafellar, \emph{Convex Analysis}.\hskip 1em plus 0.5em minus
  0.4em\relax Princeton university press, 1997.

\bibitem{hiriart2001fundamentals}
J.-B. Hiriart-Urruty and C.~Lemaréchal, \emph{Fundamentals of Convex
  Analysis}.\hskip 1em plus 0.5em minus 0.4em\relax Springer Science \&
  Business Media, 2001.

\bibitem{bauschke2011convex}
H.~H. Bauschke and P.~L. Combettes, \emph{\BIBforeignlanguage{English}{Convex
  analysis and monotone operator theory in Hilbert spaces}}.\hskip 1em plus
  0.5em minus 0.4em\relax Springer, 2011.

\bibitem{rockafellar2011variational}
R.~T. Rockafellar and R.~J.-B. Wets, \emph{Variational analysis}.\hskip 1em
  plus 0.5em minus 0.4em\relax Springer, 2011, vol. 317.

\bibitem{rockafellar1988first}
R.~T. Rockafellar, ``First- and second-order epi-differentiability in nonlinear
  programming,'' \emph{Transactions of the American Mathematical Society}, vol.
  307, no.~1, pp. 75--108, 1988.

\bibitem{rockafellar1989second}
------, ``{Second-order optimality conditions in nonlinear programming obtained
  by way of epi-derivatives},'' \emph{Mathematics of Operations Research},
  vol.~14, no.~3, pp. 462--484, 1989.

\bibitem{poliquin1992amenable}
R.~A. Poliquin and R.~T. Rockafellar, ``{Amenable functions in optimization},''
  \emph{Nonsmooth optimization: methods and applications ({E}rice, 1991)}, pp.
  338--353, 1992.

\bibitem{poliquin1995second}
------, ``Second-order nonsmooth analysis in nonlinear programming,''
  \emph{Recent advances in nonsmooth optimization}, pp. 322--349, 1995.

\bibitem{poliquin1996generalized}
------, ``Generalized {H}essian properties of regularized nonsmooth
  functions,'' \emph{SIAM Journal on Optimization}, vol.~6, no.~4, pp.
  1121--1137, 1996.

\bibitem{auslender2003asymptotic}
A.~Auslender and M.~Teboulle, \emph{Asymptotic cones and functions in
  optimization and variational inequalities}.\hskip 1em plus 0.5em minus
  0.4em\relax Springer, 2003.

\bibitem{nesterov1983method}
Y.~Nesterov, ``A method of solving a convex programming problem with
  convergence rate ${O }(1/k^2)$,'' \emph{Soviet Mathematics Doklady}, vol.~27,
  no.~2, pp. 372--376, 1983.

\bibitem{nesterov2003introductory}
------, \emph{Introductory lectures on convex optimization: A basic
  course}.\hskip 1em plus 0.5em minus 0.4em\relax Springer, 2003, vol.~87.

\bibitem{powell1970hybrid}
M.~Powell, ``A hybrid method for nonlinear equations,'' \emph{Numerical Methods
  for Nonlinear Algebraic Equations}, pp. 87--144, 1970.

\bibitem{broyden1965class}
C.~G. Broyden, ``A class of methods for solving nonlinear simultaneous
  equations,'' \emph{Mathematics of Computation}, vol.~19, no.~92, pp.
  577--593, 1965.

\bibitem{byrd1989tool}
R.~H. Byrd and J.~Nocedal, ``A tool for the analysis of quasi-{N}ewton methods
  with application to unconstrained minimization,'' \emph{SIAM Journal on
  Numerical Analysis}, vol.~26, no.~3, pp. 727--739, 1989.

\bibitem{liu1989limited}
D.~C. Liu and J.~Nocedal, ``\BIBforeignlanguage{English}{On the limited memory
  {BFGS} method for large scale optimization},''
  \emph{\BIBforeignlanguage{English}{Mathematical Programming}}, vol.~45, no.
  1-3, pp. 503--528, 1989.

\bibitem{nocedal1980updating}
J.~Nocedal, ``Updating quasi-{N}ewton matrices with limited storage,''
  \emph{Mathematics of computation}, vol.~35, no. 151, pp. 773--782, 1980.

\bibitem{nocedal2006numerical}
J.~Nocedal and S.~Wright, \emph{\BIBforeignlanguage{English}{Numerical
  Optimization}}, 2nd~ed.\hskip 1em plus 0.5em minus 0.4em\relax New York:
  Springer, 2006.

\bibitem{rockafellar1973dual}
R.~T. Rockafellar, ``A dual approach to solving nonlinear programming problems
  by unconstrained optimization,'' \emph{Mathematical Programming}, vol.~5,
  no.~1, pp. 354--373, 1973.

\bibitem{rockafellar1976augmented}
------, ``Augmented {L}agrangians and applications of the proximal point
  algorithm in convex programming,'' \emph{Mathematics of operations research},
  vol.~1, no.~2, pp. 97--116, 1976.

\bibitem{hestenes1969multiplier}
M.~R. Hestenes, ``Multiplier and gradient methods,'' \emph{Journal of
  optimization theory and applications}, vol.~4, no.~5, pp. 303--320, 1969.

\bibitem{powell1969method}
M.~J.~D. Powell, ``A method for nonlinear constraints in minimization
  problems,'' in \emph{Optimization}, R.~Fletcher, Ed.\hskip 1em plus 0.5em
  minus 0.4em\relax New York: Academic Press, 1969, pp. 283--298.

\bibitem{bertsekas2015convex}
D.~P. Bertsekas, \emph{Convex optimization algorithms}.\hskip 1em plus 0.5em
  minus 0.4em\relax Athena Scientific, 2015.

\bibitem{li1995error}
W.~Li, ``Error bounds for piecewise convex quadratic programs and
  applications,'' \emph{SIAM Journal on Control and Optimization}, vol.~33,
  no.~5, pp. 1510--1529, 1995.

\bibitem{drusvyatskiy2016error}
D.~Drusvyatskiy and A.~S. Lewis, ``Error bounds, quadratic growth, and linear
  convergence of proximal methods,'' \emph{To appear in Mathematics of
  Operations Research}, 2017.

\bibitem{dontchev2009implicit}
A.~L. Dontchev and R.~T. Rockafellar, ``Implicit functions and solution
  mappings,'' \emph{Springer Monogr. Math.}, 2009.

\bibitem{schoepfer2016linear}
F.~Sch\"opfer, ``Linear convergence of descent methods for the unconstrained
  minimization of restricted strongly convex functions,'' \emph{SIAM Journal on
  Optimization}, vol.~26, no.~3, pp. 1883--1911, 2016.

\bibitem{zhou2015unified}
Z.~Zhou and A.~M.-C. So, ``A unified approach to error bounds for structured
  convex optimization problems,'' \emph{Mathematical Programming}, vol. 165,
  no.~2, pp. 689--728, 2017.

\bibitem{aragon2008characterization}
F.~J. Arag{\'o}n~Artacho and M.~H. Geoffroy, ``Characterization of metric
  regularity of subdifferentials,'' \emph{Journal of Convex Analysis}, vol.~15,
  no.~2, pp. 365--380, 2008.

\bibitem{bernstein2009matrix}
D.~S. Bernstein, \emph{Matrix mathematics: theory, facts, and formulas}.\hskip
  1em plus 0.5em minus 0.4em\relax Princeton University Press, 2009.

\bibitem{facchinei2003finite}
F.~Facchinei and J.-S. Pang,
  \emph{\BIBforeignlanguage{English}{Finite-Dimensional Variational
  Inequalities and Complementarity Problems}}.\hskip 1em plus 0.5em minus
  0.4em\relax Springer, 2003, vol.~2.

\bibitem{themelis2016supermann}
A.~Themelis and P.~Patrinos, ``{S}uper{M}ann: a superlinearly convergent
  algorithm for finding fixed points of nonexpansive operators,'' \emph{arXiv
  preprint arXiv:1609.06955}, 2016.

\bibitem{patrinos2014accelerated}
P.~Patrinos and A.~Bemporad, ``An accelerated dual gradient-projection
  algorithm for embedded linear model predictive control,'' \emph{IEEE
  Transactions on Automatic Control}, vol.~59, no.~1, pp. 18--33, 2014.

\bibitem{giselsson2015metric}
P.~Giselsson and S.~Boyd, ``Metric selection in fast dual forward--backward
  splitting,'' \emph{Automatica}, vol.~62, pp. 1--10, 2015.

\bibitem{bemporad1997nonlinear}
A.~Bemporad, A.~Casavola, and E.~Mosca, ``Nonlinear control of constrained
  linear systems via predictive reference management,'' \emph{IEEE transactions
  on Automatic Control}, vol.~42, no.~3, pp. 340--349, 1997.

\bibitem{richter2013certification}
S.~Richter, C.~N. Jones, and M.~Morari, ``Certification aspects of the fast
  gradient method for solving the dual of parametric convex programs,''
  \emph{Mathematical Methods of Operations Research}, vol.~77, no.~3, pp.
  305--321, 2013.

\bibitem{pu2017complexity}
Y.~Pu, M.~N. Zeilinger, and C.~N. Jones, ``Complexity certification of the fast
  alternating minimization algorithm for linear {MPC},'' \emph{IEEE
  Transactions on Automatic Control}, vol.~62, no.~2, pp. 888--893, 2017.

\bibitem{ferreau2014qpoases}
H.~J. Ferreau, C.~Kirches, A.~Potschka, H.~G. Bock, and M.~Diehl, ``{qpOASES}:
  {A} parametric active-set algorithm for quadratic programming,''
  \emph{Mathematical Programming Computation}, vol.~6, no.~4, pp. 327--363,
  2014.

\bibitem{domahidi2013ecos}
A.~Domahidi, E.~Chu, and S.~Boyd, ``{ECOS}: {A}n {SOCP} solver for embedded
  systems,'' in \emph{European Control Conference (ECC)}, 2013, pp. 3071--3076.

\bibitem{toh1999sdpt3}
K.-C. Toh, M.~J. Todd, and R.~H. T{\"u}t{\"u}nc{\"u}, ``{SDPT3} -- a {MATLAB}
  software package for semidefinite programming, version 1.3,''
  \emph{Optimization methods and software}, vol.~11, no. 1-4, pp. 545--581,
  1999.

\bibitem{sturm1999using}
J.~F. Sturm, ``Using {SeDuMi} 1.02, a {MATLAB} toolbox for optimization over
  symmetric cones,'' \emph{Optimization methods and software}, vol.~11, no.
  1-4, pp. 625--653, 1999.

\bibitem{bauschke1999strong}
H.~H. Bauschke, J.~M. Borwein, and W.~Li, ``Strong conical hull intersection
  property, bounded linear regularity, {J}ameson’s property ({G}), and error
  bounds in convex optimization,'' \emph{Mathematical Programming}, vol.~86,
  no.~1, pp. 135--160, 1999.

\end{thebibliography}

	\ifieee
		\vspace{-1cm}
		\begin{IEEEbiography}[{%
			\includegraphics[width=1in,height=1.25in,clip,keepaspectratio]{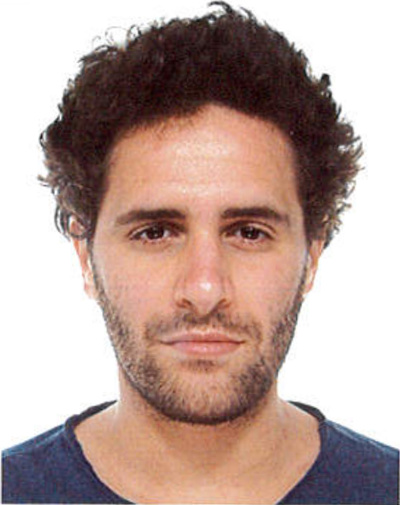}%
		}]{%
			Lorenzo Stella received the Bachelor and Master degrees in Computer Science from the University of Florence (Italy), and the Ph.D. jointly at the IMT School for Advanced Studies, Lucca (Italy) and the Department of Electrical Engineering (ESAT) of KU Leuven (Belgium).
			His research interests cover large-scale, nonsmooth optimization algorithms with applications to predictive control and machine learning problems.%
		}%
		\end{IEEEbiography}
		\vspace{-1cm}
		\begin{IEEEbiography}[{%
			\includegraphics[width=1in,height=1.25in,clip,keepaspectratio]{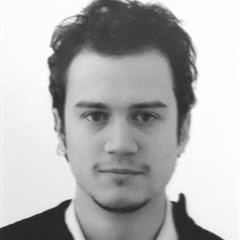}%
		}]{%
			Andreas Themelis received both Bachelor and Master degrees in Mathematics from the University of Florence, Italy, in 2010 and 2013, respectively.
			He is currently pursuing a joint Ph.D at the IMT School for Advanced Studies, Lucca (Italy) and the Department of Electrical Engineering (ESAT) of KU Leuven (Belgium).
			His research currently focuses on (non)convex nonsmooth optimization with particular interest in splitting schemes deriving from monotone operators theory, and stochastic algorithms intended for large-scale structured problems.%
		}%
		\end{IEEEbiography}
		\vspace{-1cm}
		\begin{IEEEbiography}[{%
			\includegraphics[width=1in,height=1.25in,clip,keepaspectratio]{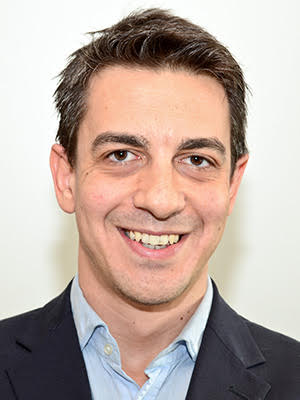}%
		}]{%
			Panagiotis (Panos) Patrinos is currently assistant professor at the Department of Electrical Engineering (ESAT) of KU Leuven, Belgium.
			He received the M.Eng. in Chemical Engineering, M.Sc. in Applied Mathematics and Ph.D. in Control and Optimization from National Technical University of Athens, Greece.
			After his Ph.D. he held postdoctoral positions at the University of Trento and IMT School of Advanced Studies Lucca, Italy, where he became an assistant professor in 2012.
			During fall/winter 2014 he held a visiting assistant professor position in the department of electrical engineering at Stanford University.
			His current research interests are in the theory and algorithms of optimization and predictive control with a focus on large-scale, distributed, stochastic and embedded optimization with a wide range of application areas including smart grids, water networks, aerospace, and machine learning.
		}
		\end{IEEEbiography}
	\fi

\end{document}